\documentclass[a4paper, 10pt, twoside]{amsart}

\usepackage{yfonts, enumitem, mathtools}
\usepackage{amsmath}
\usepackage{amssymb, latexsym, mathrsfs}
\usepackage{fancybox, graphicx, subfigure}
\usepackage{endnotes}
\usepackage[usenames, dvipsnames]{color}
\usepackage[colorlinks=true,
        raiselinks=true,
        linkcolor=MidnightBlue,
        citecolor=Mahogany,
        urlcolor=ForestGreen,
        pdfauthor=Debasish Chatterjee,
        pdftitle={Discounted stochastic shortest path problem, general state space},
        pdfkeywords={discounted stochastic shortest path, Markov control processes},
        pdfsubject={Technical Report},
        plainpages=false]{hyperref}

\usepackage[charter]{mathdesign}
\DeclareMathOperator*{\argmin}{arg\,min}

\newcommand{\R}{\ensuremath{\mathbb{R}}}
\newcommand{\N}{\ensuremath{\mathbb{N}}}
\newcommand{\Nz}{\ensuremath{\mathbb{N}_0}}
\newcommand{\posR}{\ensuremath{\R_{\ge 0}}}

\newcommand{\ra}{\ensuremath{\rightarrow}}

\newcommand{\lra}{\ensuremath{\longrightarrow}}

\newcommand{\fa}{\ensuremath{\forall\,}}
\renewcommand{\le}{\ensuremath{\leqslant}}
\renewcommand{\ge}{\ensuremath{\geqslant}}
\renewcommand{\mapsto}{\ensuremath{\longmapsto}}

\newcommand{\setmin}{\ensuremath{\!\smallsetminus}}
\newcommand{\Let}{\coloneqq}

\newcommand{\norm}[1]{\ensuremath{\left\lVert #1 \right\rVert}}
\newcommand{\Borelsigalg}[1]{\ensuremath{\mathfrak{B}\!\left(#1\right)}}

\newcommand{\indic}[1]{\ensuremath{\boldsymbol{1}_{#1}}}
\newcommand{\abs}[1]{\ensuremath{\left\lvert{#1}\right\rvert}}

\newcommand{\mrm}[1]{\ensuremath{\mathrm{#1}}}
\newcommand{\mc}[1]{\ensuremath{\mathcal{#1}}}

\newcommand{\Lp}[1]{\ensuremath{\boldsymbol{L}_{#1}}}

\newcommand{\secref}[1]{\S\ref{#1}}

\newcommand{\EE}{\ensuremath{\mathsf{E}}}
\newcommand{\PP}{\ensuremath{\mathsf{P}}}

\newcommand{\sigalg}{\ensuremath{\mathfrak{F}}}

\newcommand{\xz}{\ensuremath{x_0}}

\newcommand{\ol}{\overline}

\newcommand{\wt}{\widetilde}
\newcommand{\wh}{\widehat}

\renewcommand{\subset}{\ensuremath{\subseteq}}

\newcommand{\mn}{\ensuremath{\wedge}}

\newcommand{\RemarkEnd}{\hspace{\stretch{1}}{$\vartriangleleft$}}
\newcommand{\ExampleEnd}{\hspace{\stretch{1}}{$\triangle$}}
\newcommand{\DefEnd}{\hspace{\stretch{1}}{$\Diamond$}}
\newcommand{\AssumptionEnd}{\hspace{\stretch{1}}{$\diamondsuit$}}
\newcommand{\vphi}{\ensuremath{\varphi}}

\numberwithin{equation}{section}
\swapnumbers
\newtheorem{theorem}[equation]{Theorem}
\newtheorem{corollary}[equation]{Corollary}
\newtheorem{lemma}[equation]{Lemma}
\newtheorem{proposition}[equation]{Proposition}

\theoremstyle{definition}
\newtheorem{defn}[equation]{Definition}
\theoremstyle{remark}
\newtheorem{remark}[equation]{Remark}
\newtheorem{example}[equation]{Example}
\newtheorem{prgr}[equation]{}
\newtheorem{assumption}[equation]{Assumption}

\allowdisplaybreaks

\title[On Stochastic Control up to a Hitting Time]{Stochastic Control up to a Hitting Time: Optimality and Rolling-horizon Implementation}
\author[D.~Chatterjee]{Debasish Chatterjee}
\author[E.~Cinquemani]{Eugenio Cinquemani}
\author[G.~Chaloulos]{Georgios Chaloulos}
\author[J.~Lygeros]{John Lygeros}
\address{Automatic Control Laboratory, Physikstrasse 3, ETH Z\"urich, 8092 Z\"urich, Switzerland}
\email{\{chatterjee,cinquemani,chaloulos,lygeros\}@control.ee.ethz.ch}
\urladdr{\url{http://control.ee.ethz.ch}}

\date{\today}

\begin{document}

	\begin{abstract}
		We present a dynamic programming-based solution to a stochastic optimal control problem up to a hitting time for a discrete-time Markov control process. First we determine an optimal control policy to steer the process toward a compact target set while simultaneously minimizing an expected discounted cost. We then provide a rolling-horizon strategy for approximating the optimal policy, together with quantitative characterization of its sub-optimality with respect to the optimal policy. Finally we address related issues of asymptotic discount-optimality of the value-iteration policy.
	\end{abstract}

	\maketitle

	\section{Introduction}
	\label{s:intro}
		Optimal control of Markov control processes (MCP) up to an exit time is a problem with a long and rich history. It has mostly been studied as the minimization of an expected undiscounted cost until the first time that the state enters a given target set, see e.g.,~\cite[Chapter~II]{ref:borkarTopicsControlledMC},~\cite[Chapter~8]{ref:hernandez-lerma2}, and the references therein. In particular, if a unit cost is incurred as long as the state is outside the target set, then the problem of minimizing the cost accumulated until the state enters the target is known variously as the \textsl{pursuit problem}~\cite{ref:eatonzadeh62}, \textsl{transient programming}~\cite{ref:whittleOptimization}, the \textsl{first passage problem}~\cite{ref:dermanMDP, ref:kushnerIntroStochControl}, the \textsl{stochastic shortest path problem}~\cite{ref:bertsekasDP2}, and \textsl{control up to an exit time}~\cite{ref:borkarConvexAnalyticApproach, ref:borkarTopicsControlledMC, ref:kestenMCP}. These articles deal with at most countable state and action spaces. The problem of optimally controlling a system until an exit time from a given set has gained significance in financial and insurance mathematics, see, e.g., \cite{ref:boda04, ref:schmidliInsurance}.

		Our interest in this problem stems from our attempts to develop a general theory of stochastic model-predictive control (MPC). In its bare essentials, deterministic MPC~\cite{ref:maciejowskibk} consists of two steps: (i) solving a finite-horizon optimal control problem with constraints on the state and the controlled inputs to get an optimal policy, and (ii) applying a controller derived from the policy obtained in step (i) in a rolling-horizon fashion. Theoretical foundation of stochastic MPC is still in its infancy, see~\cite{ref:PrimbsSung09, ref:bertsimas2007, ref:vanHessem2006, ref:kouvaritakissMPCIneqconstraints, ref:batinaPhDthesis} and the references therein for some related work. In view of its close relationship with applications, any satisfactory theory of stochastic MPC must necessarily take into account its practical aspects. In this context an examination of a standard linear system with constrained controlled inputs affected by independent and identically distributed (i.i.d.)\ unbounded (e.g., Gaussian) disturbance inputs shows that no control policy can ensure that with probability one the state stays confined to a bounded \emph{safe set} for all instants of time. This is because the noise is unbounded and the samples are independent of each other. Although disturbances are not likely to be unbounded in practice, assigning an a priori bound seems to demand considerable insight. In case a bounded-noise model is adopted, existing robust MPC techniques~\cite{ref:bemporad1999rmp, ref:blanchini1999sic} may be applied, in which the central idea is to synthesize a controller based on the bounds of the noise such that the target set becomes invariant with respect to the closed-loop dynamics. However, since the optimal policy is based on a worst-case analysis, it usually leads to rather conservative controllers and sometimes even to infeasibility. Moreover, complexity of the optimization problem grows rapidly (typically exponentially) with the optimization horizon. An alternative is to replace the hard constraints by probabilistic (soft) ones. The idea is to find a policy that guarantees that the state constraints are satisfied with high probability over a sufficiently long time horizon. While this approach may improve feasibility aspects of the problem, it does not address the issue of what actions should be taken once the state violates the constraints. See~\cite{ref:hokayemcdc09, ref:smpcbnddu, ref:accl08} for recent results in this direction.

		In view of the above considerations, developing recovery strategies appears to be a necessary step. Such a strategy is to be activated once the state violates the constraints and to be deactivated whenever the system returns to the safe set. In general, a recovery strategy must drive the system quickly to the safe set while simultaneously meeting other performance objectives. In the context of MPC, two merits are immediate: (a) once the constraints are transgressed, appropriate actions can be taken to bring the state back to the safe set quickly and optimally, and (b) if the original problem is posed with hard constraints on the state, in view of (a) they may be relaxed to probabilistic ones to improve feasibility.

		In this article we address the problem of synthesizing optimal recovery strategies. We formulate the problem as the minimization of an expected discounted cost until the state enters the safe set. An almost customary assumption in the literature (see, e.g.,~\cite{ref:hindererAbsorbingSet} and the references therein,) concerned with stochastic optimal control up to an exit time is that the target set is absorbing. That is, there exists a control policy that makes the target set invariant with respect to the closed-loop stochastic dynamics. This is rather restrictive for MPC problems---it is invalid, for instance, in the very simple case of a linear controlled system with i.i.d.\ Gaussian noise inputs. We do not make this assumption, for, as mentioned above, our primary motivation for solving this problem is precisely to deal with the case that the target set is not absorbing. As a result of this, it turns out that the dynamic programming equations involve integration over subsets of the state-space and therefore are difficult to solve. At present there is no established method to solve such equations in uncountable state-spaces. However, in finite state-space cases tractable approximate dynamic programming methods~\cite{ref:bertsekasNDP, ref:powellADP} may be employed to arrive at suboptimal but efficient policies.

		This article unfolds as follows. In~\secref{s:prelims} we define the general setting of the problem, namely, Markov control processes on Polish spaces, their transition kernels and the main types of control strategies. In~\secref{s:EDC} we establish our main Theorem~\ref{t:EDC} under standard mild hypotheses. This result guarantees the existence of a deterministic stationary policy that leads to the minimal cost and also provides a Bellman equation that the value function must satisfy. A contraction mapping approach to the problem is pursued in~\secref{s:contr} under the (standard) assumption that the cost-per-stage function satisfies certain growth-rate conditions. The main result (Proposition~\ref{p:Tfp}) of this section asserts both the existence and uniqueness of the optimal value function. Asymptotic discount-optimality of the value-iteration policy is investigated in~\secref{s:ado} under two different sets of hypotheses; in particular, the results of this section show that rolling-horizon strategy approaches optimality as the length of the horizon window increases to infinity. A rolling-horizon strategy corresponding to our optimal control problem is developed in~\secref{s:rh}; in Theorem~\ref{t:rh} we establish quantitative bounds on the degree of sub-optimality of the rolling-horizon strategy with respect to the optimal policy. We conclude in~\secref{s:concl} with a discussion of future work. The state and control/action sets are assumed to be Borel subsets of Polish spaces.

	\section{Preliminaries}
	\label{s:prelims}
		We employ the following standard notations. Let $\N$ denote the natural numbers $\{1, 2, \ldots\}$, and $\Nz$ denote the nonnegative integers $\{0\}\cup\N$. Let $\indic{A}(\cdot)$ be the standard indicator function of a set $A$, i.e., $\indic{A}(\xi) = 1$ if $\xi\in A$ and $0$ otherwise. For two real numbers $a$ and $b$, let $a\mn b \Let \min\{a, b\}$.

		Given a nonempty Borel set $X$ (i.e., a Borel subset of a Polish space), its Borel $\sigma$-algebra is denoted by $\Borelsigalg{X}$. By convention ``measurable'' means ``Borel-measurable'' in the sequel. If $X$ and $Y$ are nonempty Borel spaces, a \emph{stochastic kernel} on $X$ given $Y$ is a mapping $Q(\cdot|\cdot)$ such that $Q(\cdot|y)$ is a probability measure on $X$ for each fixed $y\in Y$, and $Q(B|\cdot)$ is a measurable function on $Y$ for each fixed $B\in\Borelsigalg X$. We let $\mc P(X|Y)$ be the family of all stochastic kernels on $X$ given $Y$.% We say that $Q(\cdot|\cdot)$ is a \emph{sub-stochastic kernel} on $X$ given $Y$ if $Q(B|\cdot)$ is a measurable function on $Y$ for each $B\in\Borelsigalg{X}$, and $Q(\cdot|y)$ is a measure on $X$ with $Q(X|y) \le 1$ for each $y\in Y$.

		We briefly recall some standard definitions.

		\begin{defn}
		\label{d:mcm}
			A \emph{Markov control model} is a five-tuple
			\begin{equation}
				\label{e:mmodel}
				\bigl(X, A, \{A(x)\mid x\in X\}, Q, c\bigr)
			\end{equation}
			consisting of a nonempty Borel space $X$ called the \emph{state space}, a nonempty Borel space $A$ called the \emph{control} or \emph{action set}, a family $\{A(x)\mid x\in X\}$ of nonempty measurable subsets $A(x)$ of $A$, where $A(x)$ denotes the set of \emph{feasible controls} or \emph{actions} when the system is in state $x\in X$, and with the property that the set $\mathbb K \Let \bigl\{(x, a)\big|x\in X, a\in A(x)\bigr\}$ of feasible state-action pairs is a measurable subset of $X\times A$, a stochastic kernel $Q$ on $X$ given $\mathbb K$ called the \emph{transition law}, and a measurable function $c:\mathbb K\lra \R$ called the \emph{cost-per-stage function}.\DefEnd
		\end{defn}

		\begin{assumption}
		\label{a:basic}
			The set $\mathbb K$ of feasible state-action pairs contains the graph of a measurable function from $X$ to $A$.\AssumptionEnd%That is, there exists a measurable function $f:X\lra A$ such that $f(x)\in A(x)$ for all $x\in X$.
		\end{assumption}

		We let $\Pi$, $\Pi_{RM}$, $\Pi_{DM}$ and $\Pi_{DS}$ denote the set of all randomized and history-dependent admissible policies, randomized Markov, deterministic Markov and deterministic stationary policies, respectively. For further details and notations on policies see, e.g.,~\cite{ref:hernandez-lerma1}. Consider the Markov control model~\eqref{e:mmodel}, and for each $i=0, 1, \ldots,$ define the space $H_i$ of \emph{admissible histories} up to time $i$ as $H_0 \Let X$, and $H_i \Let \mathbb K^i\times X = \mathbb K\times H_{i-1}$ for $i\in \N$. A generic element $h_i$ of $H_i$, called an admissible $i$-history is a vector of the form $h_i = (x_0, a_0, \ldots, x_{i-1}, a_{i-1}, x_i)$, with $(x_j, a_j)\in\mathbb K$ for $j=0, \ldots, i-1$ and $x_i\in X$. Hereafter we let the $\sigma$-algebra generated by the history $h_i$ be denoted by $\sigalg_i$, $i\in\Nz$. %A policy $\pi = (\pi_i)_{i\in\Nz}$ is a sequence $(a_i)_{i\in\Nz}$ of $A$-valued random variables, called actions or controls, such that for every $i$-history $h_i, i\in\Nz$, the law of $a_i$ is $\pi_i(\cdot|h_i)$, with support of $\pi_i(\cdot|h_i)$ contained in $A(x_i)$, the set of feasible actions in the state $x_i$. 
		Let $(\Omega, \sigalg)$ be the measurable space consisting of the (canonical) sample space $\Omega \Let \ol H_\infty = (X\times A)^\infty$, and $\sigalg$ is the corresponding product $\sigma$-algebra. %For $\omega = (x_0, a_0, x_1, a_1, \ldots)\in\Omega$, the projections $x_i$ and $a_i$ from $\Omega$ to the sets $X$ and $A$ are called \emph{state} and \emph{control} (or \emph{action}) variables, respectively. 
		Let $\pi = (\pi_i)_{i\in\Nz}$ be an arbitrary control policy and $\nu$ an arbitrary probability measure on $X$, referred to as the initial distribution. By a theorem of Ionescu-Tulcea~\cite[Chapter 3, \S4, Theorem~5]{ref:raoProbTheo}, there exists a unique probability measure $\mathsf P_\nu^\pi$ on $(\Omega, \sigalg)$ supported on $H^\infty$, and such that for all $B\in\Borelsigalg X$, $C\in\Borelsigalg A$, and $h_i\in H_i$, $i\in\Nz$, $\mathsf P_\nu^\pi\bigl(\xz\in B\bigr) = \nu(B)$ and
		\begin{subequations}
		%\label{e:probmeasure}
		\begin{align}
			\label{e:actiontrans}
			\mathsf P_\nu^\pi\bigl(a_i\in C\,\big|\, h_i\bigr) & = \pi_i\bigl(C\,\big|\, h_i\bigr)\\
			\label{e:statetrans}
			\mathsf P_\nu^\pi\bigl(x_{i+1}\in B\,\big|\, h_i, a_i\bigr) & = Q\bigl(B\,\big|\, x_i, a_i\bigr).
		\end{align}
		\end{subequations}
		The stochastic process $\bigl(\Omega, \sigalg, \mathsf P_\nu^\pi, (x_i)_{i\in\Nz}\bigr)$ is called a discrete-time \emph{Markov control process}. Let $\Phi$ denote the set of stochastic kernels $\vphi$ in $\mathcal P(A| X)$ such that $\vphi(A(x)| x) = 1$ for all $x\in X$, and let $\mathbb F$ denote the set of all measurable functions $f:X\lra A$ satisfying $f(x)\in A(x)$ for all $x\in X$. The functions in $\mathbb F$ are called \emph{selectors} of the set-valued mapping $X\ni x\mapsto A(x)\subset A$.
		%Recall that a policy $\pi = (\pi_i)_{i\in\Nz}\in\Pi$ is said to be \emph{randomized Markov} if there exists a sequence $(\vphi_i)_{i\in\Nz}$ of stochastic kernels $\vphi_i\in\Phi$ such that $\pi_i(\cdot| h_i) = \vphi_i(\cdot| x_i)\;\; \fa h_i\in H_i, \;i\in\Nz$, \emph{deterministic Markov} if there exists a sequence $(f_i)_{i\in\Nz}$ of functions $f_i\in\mathbb F$ such that $\pi_i(\cdot| h_i) = \delta_{f(x_i)}(\cdot)$, and \emph{deterministic stationary} if there exists a function $f\in\mathbb F$ such that $\pi_i(\cdot| h_i) = \delta_{f(x_i)}(\cdot)$.

		%We note that the process $\bigl(\Omega, \sigalg, \mathsf P_\nu^\pi, (x_i)_{i\in\Nz}\bigr)$ is not necessarily Markovian due to the dependence on the entire history $h_i$ in~\eqref{e:actiontrans}; however, if $(\pi_i)_{i\in\Nz}$ is restricted to randomized Markov policies, then $(x_i)_{i\in\Nz}$ is a Markov process, as established in~\cite[Proposition~2.3.5]{ref:hernandez-lerma1}.

		The transition kernel $Q$ in~\eqref{e:statetrans} under a policy $\pi \Let (\vphi_i)_{i\in\Nz}\in\Pi_{RM}$ is given by $\bigl(Q(\cdot|\cdot, \vphi_i)\bigr)_{i\in\Nz}$, defined as $\Borelsigalg{X}\times X\ni (B, x)\mapsto Q(B|x, \vphi_i(x)) \Let \int_{A(x)}\vphi_i(\mrm da|x) Q(B|x, a)$. Occasionally we suppress the dependence of $\vphi_i$ on $x$ and write $Q(B|x, \vphi_i)$ in place of $Q(B|x, \vphi_i(x))$. The cost-per-stage function at the $j$-th stage under a policy $(\vphi_i)_{i\in\Nz}$ is written as $c(x_j, \vphi_j) \Let \int_{A(x_j)} \vphi_j(\mrm da|x_j)c(x_j, a)$. We simply write $\vphi^\infty$ and $f^\infty$, respectively, for policies $(\vphi, \vphi, \ldots)\in\Pi_{RS}$ and $(f, f, \ldots)\in \Pi_{DS}$.

		%\begin{proposition}[{\cite[Proposition~2.3.5]{ref:hernandez-lerma1}}]
		%	\label{p:Markovprop}
		%	Let $\nu$ be an arbitrary initial distribution. If $\pi = (\vphi_i)_{i\in\Nz}$ is a randomized or deterministic Markov policy, then $(x_i)_{i\in\Nz}$ is a inhomogeneous Markov process with transition kernel $(Q(\cdot| \cdot, \vphi_i)$ at the $i$-th step. In particular, if $\pi = \vphi^\infty$ and $\pi = f^\infty$ are a stationary randomized and a deterministic stationary policy, respectively, then $(x_i)_{i\in\Nz}$ is a time-homogeneous Markov process with corresponding transition kernels $Q(\cdot| \cdot, \vphi)$ and $Q(\cdot| \cdot, f)$ at each step, respectively.
		%\end{proposition}

		Since we shall be exclusively concerned with Markov policies and its subclasses, in the sequel we use the notation $\Pi$ for the class of all randomized Markov strategies.

	\section{Expected Discounted Cost up to the first Exit Time}
	\label{s:EDC}
		Let $K\subset X$ be a measurable set, $x_0 = x\in X$ and let $\tau \Let \inf\bigl\{i\in\Nz\big| x_i\in K\bigr\}$.\footnote{As usual the infimum over an empty set is taken to be $+\infty$.} We note that $\tau$ is an $(\sigalg_i)_{i\in\Nz}$-stopping time. Let us define
		\[
			V(\pi, x) \Let \mathsf E_x^\pi\!\left[\sum_{i=0}^{\tau-1} \alpha^i c(x_i, a_i)\right], \qquad \alpha\in\:]0, 1[,
		\]
		as the \emph{$\alpha$-discounted expected cost} under policy $\pi\in\Pi$ corresponding to the Markov control process $\bigl(\Omega, \sigalg, \mathsf P_\nu^\pi, (x_i)_{i\in\Nz}\bigr)$.\footnote{We employ the standard convention that a summation from a higher to a lower index is defined to be $0$.} Our objective is to minimize $V(\pi,x)$ over a class of control policies $\Pi$, i.e., find the $\alpha$-discount value function
		\begin{align}
			\label{e:problem}
			V^\star(x) \Let \inf_{\pi\in\Pi} V(\pi, x) = \inf_{\pi\in\Pi}\mathsf E_x^\pi\!\left[\sum_{i=0}^{\tau-1} \alpha^i c(x_i, a_i)\right], \qquad \alpha\in\:]0, 1[.
		\end{align}
		A policy that attains the infimum above is said to be \emph{$\alpha$-discount optimal}.
		
		\begin{remark}
			As mentioned in the introduction, the optimization problem~\eqref{e:problem} with $\alpha = 1$ and the cost-per-stage function $c(x, a) = \indic{X\setmin K}(x)$ is known as the stochastic shortest path problem. The objective of this problem is to drive the state to a desired set ($K$ in our case) as soon as possible, and the expected cost $V_{\text{ssp}}(\pi, x)$ for a policy $\pi$ corresponding to the above cost-per-stage function is readily seen to be $\mathsf E^\pi_x\bigl[\tau\bigr]$. In this light we observe that the minimization problem in~\eqref{e:problem} with the cost-per-stage function $c(x, a) = \indic{X\setmin K}(x)$ can be viewed as a discounted stochastic shortest path problem. It follows immediately that the corresponding expected cost $V_{\text{dssp}}(\pi, x)$ is $\bigl(1-\mathsf E^\pi_x\bigl[\alpha^\tau\bigr]\bigr)/(1-\alpha)$. Note that the minimization of $V_{\text{dssp}}(\pi, x)$ over a class of policies is always well-defined for $\alpha < 1$. Moreover, because of the monotonic behavior of the map $]0, 1[\;\ni\alpha\mapsto \bigl(1-\mathsf E^\pi_x\bigl[\alpha^\tau\bigr]\bigr)/(1-\alpha)$, one may hope to get a good approximation of the original stochastic shortest path problem. However, pathological examples can be constructed to show that a solution to the stochastic shortest path problem may not exist, whereas minimization of $V_{\text{dssp}}(\pi, x)$ is always well defined, although in either case the state may never reach the desired set $K$ almost surely-$\PP^\pi_x$.\RemarkEnd
		\end{remark}

		\begin{remark}
			\label{r:diffc}
			%If we take the particular case of $c(x, a) = \indic{X\setmin K}(x)$, then we obtain the problem of minimizing the time taken by the process $(x_i)_{i\in\Nz}$ to hit the set $K$ for the first time, which is the stochastic shortest time to hit $K$ problem. Also 
			Given a cost-per-stage function $c$ on $\mathbb K$, one can redefine it to be $c'(x, a) \Let c(x, a)\indic{X\setmin K}(x)$ to turn the problem~\eqref{e:problem} into the minimization of $\mathsf E_x^\pi\!\left[\sum_{i=0}^\tau \alpha^i c'(x_i, a_i)\right]$ for $\alpha\in\:]0, 1[$. This cost functional can be equivalently written as an infinite horizon cost functional, as in $\mathsf E_x^\pi\!\left[\sum_{i=0}^\infty \alpha^i c'(x_i, a_i)\indic{\{i \le \tau\}}\right]$, or as in $\mathsf E_x^\pi\!\left[\sum_{i=0}^\infty \alpha^i c(x_i, a_i)\indic{\{i < \tau\}}\right]$. However, the absence of a policy that guarantees that $(x_i)_{i\in\Nz}$ stays inside $K$ for all time after $\tau$ necessarily means that the problem~\eqref{e:problem} corresponding to the Markov control model in Definition~\ref{d:mcm} is not equivalent to the minimization of the infinite horizon cost functional $\mathsf E_x^\pi\!\left[\sum_{i=0}^\infty \alpha^i c'(x_i, a_i)\right]$.\RemarkEnd
		\end{remark}

		\begin{prgr}
		\label{pgr:policies}
			\emph{A word about admissible policies.} It is clear at once that the class of admissible policies for the problem~\eqref{e:problem} is different from the classes considered in~\secref{s:prelims}. Indeed, since the process is killed at the stopping time $\tau$, it follows that the class of admissible policies should also be truncated at the stage $\tau-1$. For a given stage $t\in\Nz$ we define the $t$-th policy element $\pi_t$ only on the set $\{t < \tau\}$. Note that with this definition $\pi_t$ becomes a $\sigalg_{t\mn\tau}$-measurable randomized control (in general). It is also immediate from the definition of $\tau$ that if the initial condition $x$ is inside $K$, then the set of admissible policies is empty; indeed, in this case $\tau = 0$, and there is no control needed. In other words, the domain of $\pi_t$ is contained in the ``spatial'' region $\bigl\{(x, a)\in\mathbb K\,\big|\,x\in X\setmin K, a\in A(x)\bigr\}$; since $\pi_t$ is not defined on $K$, this is equivalent to $\pi_t$ being well-defined on $\{t < \tau\}$.
		\end{prgr}

		\begin{prgr}
		\label{pgr:convention}
			\emph{Some re-definitions.} To simplify the formulas from now on we let the cost-per-stage function to be defined on $X\setmin K$. With this convention in place our problem~\eqref{e:problem} can be posed as the minimization of $\EE^\pi_x\bigl[\sum_{i=0}^{\tau-1} \alpha^i c(x_{i}, a_{i})\bigr]$ over admissible policies. Also, henceforth we redefine the set $\mathbb K$ of state-action pairs to be $\mathbb K \Let \bigl\{(x, a)\in X\times A\,\big|\,x\in X\setmin K, a\in A(x)\bigr\}$, and we note that this new set is a measurable subset of the original set of state-action pairs. Also, we let $\mathbb F$ be the set of selectors of the set-valued mapping $X\setmin K\ni x\mapsto A(x)\subset A$.
		\end{prgr}

		%\begin{prgr}
		%	\label{prgr:fndef}
		Recall that a function $g:\mathbb K\lra \R$ is said to be \emph{inf-compact on $\mathbb K$} if for every $x\in X$ and $r\in \R$ the set $\bigl\{a\in A(x)\big| g(x, a) \le r\bigr\}$ is compact. A transition kernel $Q$ on a measurable space $X$ given another measurable space $Y$ is said to be \emph{strongly Feller} (or \emph{strongly continuous}) if the mapping $y\mapsto \int_X g(x) Q(\mrm dx| y)$ is continuous and bounded for every measurable and bounded function $g:X\lra\R$. A function $g:\mathbb K\lra\R$ is \emph{lower semicontinuous} (l.s.c.) if for every sequence $(x_j, a_j)_{j\in\N}\subset\mathbb K$ converging to $(x, a)\in\mathbb K$, we have $\liminf_{j\ra\infty} g(x_j, a_j) \ge g(x, a)$; or, equivalently, if for every $r\in\R$, the set $\bigl\{(x, a)\in\mathbb K\big| g(x, a) \le r\bigr\}$ is closed in $\mathbb K$.
		%\end{prgr}

		\begin{assumption}
			\label{a:key}
			In addition to Assumption~\ref{a:basic}, we stipulate that
			\begin{enumerate}[align=right, leftmargin=*, widest=iii, label=(\roman*)]
				\item the set $A(x)$ is compact for every $x\in X$,
				\item the cost-per-stage $c$ is lower semicontinuous, nonnegative, and inf-compact on $\mathbb K$, and
				\item the transition kernel $Q$ is strongly Feller.\AssumptionEnd
			\end{enumerate}
		\end{assumption}

		The following is our main result on expected discounted cost up to the first time $\tau$ to hit $K$; a proof is presented later in this section.

		\begin{theorem}
			\label{t:EDC}
			Suppose that Assumption {\rm \ref{a:key}} holds. Then
			\begin{enumerate}[label=\emph{(\roman*)}, align=right, leftmargin=*, widest=iii]
				\item The $\alpha$-discount value function $V^\star$ is the (positive) minimal measurable solution to the $\alpha$-discounted cost optimality equation ($\alpha$-DCOE)
				\begin{equation}
					\label{e:alphadcoe}
					\xi(x) = \min_{A(x)}\left[c(x, a) + \alpha\int_{X\setmin K} Q(\mrm dy| x, a)\:\xi(y)\right]\qquad \fa x\in X\setmin K.
				\end{equation}
				\item There exists a selector $f_\star\in\mathbb F$ such that $f_\star(x)\in A(x)$, $x\in X\setmin K$, attains the minimum in~\eqref{e:alphadcoe}, i.e.,
				\begin{equation}
					\label{e:alphado}
					V^\star(x) = c(x, f_\star) + \alpha\int_{X\setmin K}Q(\mrm dy| x, f_\star)\:V^\star(y)\qquad \fa x\in X\setmin K,
				\end{equation}
				and the deterministic stationary policy $f_\star^\infty$ is $\alpha$-discount optimal; conversely, if $f_\star^\infty\in\Pi_{DS}$ is $\alpha$-discount optimal, then it satisfies~\eqref{e:alphado}.
				%\item If an $\alpha$-discount optimal policy exists, then there exists one that is deterministic stationary.
			\end{enumerate}
		\end{theorem}

		We observe that Theorem~\ref{t:EDC} does not assert that the optimal value function $V^\star$ is unique in any sense. In~\secref{s:contr} we prove a result (Proposition~\ref{p:Tfp}) under additional hypotheses that guarantees uniqueness of $V^\star$.

		Since we do not assume that the cost-per-stage function $c$ is bounded, a useful approach is to consider the $\alpha$-\emph{value iteration} ($\alpha$-VI) \emph{functions} defined by
		\begin{equation}
		\label{e:VI}
		\begin{cases}
			v_0(x) = 0,\\
			v_n(x) = \displaystyle{\min_{A(x)}\left[c(x, a) + \alpha \int_{X\setmin K} Q(\mrm dy|x, a)\: v_{n-1}(y)\right]},
		\end{cases}
		n\in\N,\;\; x\in X\setmin K.
		\end{equation}
		Of course we have to demonstrate that $V^\star(x) = \lim_{n\ra\infty} v_n(x)$ for all $x\in X$.

		The functions $v_n$, $n\in\N$, may be identified with the optimal cost function for the minimization of the process stopped at the $n\mn(\tau-1)$-th step, i.e.,
		\[
			v_n(x) = \inf_{\pi\in\Pi} \mathsf E_x^\pi\!\left[\sum_{i=0}^{(n-1)\mn(\tau-1)} \alpha^i c(x_i, a_i)\right].
		\]
		To get an intuitive idea, fix a deterministic Markov policy $\pi = (\pi_i)_{i\in\Nz}$, and take the first iterate $v_1$. From~\eqref{e:VI} it is immediately clear that $v_1(x) = \min_{a\in A(x)} c(x, a)$ if $x\not\in K$, and not defined otherwise. For the second iterate, we have
		\begin{align*}
			v_2(x) & = \inf_{\pi\in\Pi} \mathsf E_x^\pi\!\left[\sum_{i=0}^{1\mn (\tau-1)} \alpha^i c(x_i, a_i)\right]\\
			& = \inf_{\pi\in\Pi}\left(c(x, \pi_0(x)) + \alpha\!\int_X \! Q(\mrm d\xi_1|x, \pi_1(x))\indic{X\setmin K}(\xi_1) c(\xi_1, \pi_1)\right).
		\end{align*}
		Note that only those sample paths that do not enter $K$ at the first step contribute to the cost at the second stage. This property is ensured by the indicator function that appears on the right-hand side of the last equality above.

		\begin{example}
			Let $(x_i)_{i\in\Nz}$ be a Markov chain with state-space $X = \{1, 2, \ldots, m\}$ and transition probability matrix $Q = [q_{ij}(a)]_{m\times m}$, where the argument of $q_{ij}$ depicts the dependence on the action $a\in A$ with $A$ being a compact subset of $\R$. Let $K = \{1, 2, \ldots, m'\}$ for $m' < m$, fix $\alpha\in\;]0, 1[$ and let $c(x, a) \Let \indic{X\setmin K}(x)$. Suppose further that $\inf_{a}q_{ij}(a) > 0$ for all $i, j\in X$; this means, in particular, that the target set $K$ cannot be absorbing for any deterministic stationary policy. Our objective is to find an optimal policy corresponding to the the minimal cost~\eqref{e:alphado}. The optimal value function $V^\star$ is $0$ on $K$ and for every $i\in \{m'+1, \ldots, m\}$ we have $V^\star(i) = \min_{a\in A(i)}\bigl[\indic{X\setmin K}(i) + \alpha\int_{X\setmin K}Q(\mrm dy|i, a) V^\star(y)\bigr] = 1 + \alpha\min_{a\in A(i)}\sum_{j=m'+1}^m q_{ij}(a) V^\star(j)$. The most elementary case is that of $m' = m - 1$; then $V^\star(m) = 1 + \alpha\min_{a\in A(m)}q_{mm}(a)V^\star(m)$, and given a sufficiently regular function $q_{mm}(\cdot)$ this can be solved at once to get $V^\star(m)$, which characterizes the function (vector) $V^\star$ completely. The optimal policy in this case is $f(m) \in \argmin_{a\in A(m)} q_{mm}(a)V^\star(m)$; if the function $q_{mm}(\cdot)$ is convex, then the minimum is attained on $A$ and thus leads to a unique optimal policy.\ExampleEnd
		\end{example}

		\subsection*{Proof of Theorem~\ref{t:EDC}}
		%\label{s:proofsEDC}
		Recall from paragraph~\ref{pgr:convention} that $c$ is defined on $X\setmin K$, $\mathbb K = \bigl\{(x, a)\in X\times A\,\big|\, x\in X\setmin K, a\in A(x)\bigr\}$ and $\mathbb F$ is the set of selectors of the set-valued map $X\setmin K\ni x\mapsto A(x)\subset A$. We begin with a sequence of Lemmas.

		\begin{lemma}[{\cite[Lemma~4.2.4]{ref:hernandez-lerma1}}]
			\label{l:keyconvergence}
			Let the functions $u:\mathbb K\lra\R$ and $u_i:\mathbb K\lra\R$, $i\in\N$, be l.s.c., inf-compact and bounded below. If $u_i\uparrow u$, then
			\[
				\lim_{i\ra\infty}\min_{A(x)} u_i(x, a) = \min_{A(x)} u(x, a)\qquad \fa x\in X.
			\]
		\end{lemma}

		\begin{lemma}[Adapted from \cite{ref:riederselectors}]
			\label{l:basicselector}
			Suppose that
			\begin{itemize}[label=\textbullet, align=right, leftmargin=*]
				\item $A(x)$ is compact for each $x\in X\setmin K$ and $\mathbb K$ is a measurable subset of $(X\setmin K)\times A$, and
				\item $v:\mathbb K\lra\posR$ is a measurable inf-compact function, $v(x, \cdot)$ is l.s.c.\ on $A(x)$ for each $x\in X$.
			\end{itemize}
			Then there exists a selector $f_\star\in\mathbb F$ such that
			\[
				v(x, f_\star(x)) = v^\star(x) \Let \min_{A(x)}v(x, a)\qquad \fa x\in X\setmin K,
			\]
			and $v^\star$ is a measurable function.
		\end{lemma}

		\begin{defn}
			Let $\Lp 0(X\setmin K)^+$ denote the convex cone of nonnegative extended real-valued measurable functions on $X\setmin K$, and for every $u\in \Lp 0(X\setmin K)^+$ let us define the map $T u$ by
			\begin{equation}
				\label{e:Tdef}
				X\setmin K\ni x\mapsto T u(x) \Let \inf_{A(x)}\left[c(x, a) + \alpha\int_{X\setmin K}Q(\mrm dy|x, a)\:u(y)\right].
			\end{equation}
			The map $T$ is the \emph{dynamic programming operator} corresponding to our problem~\eqref{e:problem}.\DefEnd
		\end{defn}

		Having defined the dynamic programming operator $T$ above, it is important to distinguish conditions under which the function $T u$ is measurable for $u\in\Lp 0(X\setmin K)^+$. We have the following lemma.

		\begin{lemma}
			\label{l:selector}
			Under Assumption {\rm \ref{a:key}}, the mapping $T$ in~\eqref{e:Tdef} takes $\Lp 0(X\setmin K)^+$ into itself. Moreover, there exists a selector $f\in \mathbb F$ such that $T u$ defined in~\eqref{e:Tdef} satisfies
			\begin{equation}
				\label{e:Tviaselector}
				T u(x) = c(x, f) + \alpha\int_{X\setmin K} Q(\mrm dy|x, f)\:u(y)\qquad \fa x\in X\setmin K.
			\end{equation}
		\end{lemma}
		\begin{proof}
			Fix $u\in\Lp 0(X\setmin K)^+$. The strong-Feller property of $Q$ on $\mathbb K$ and lower-semicontinuity of the cost-per-stage function $c$ defined on $K$ show that the map 
			\[
				\mathbb K\ni (x, a)\mapsto T'u(x, a) \Let c(x, a) + \alpha\int_{X}Q(\mrm dy|x, a)\;\indic{X\setmin K}(y) u(y)
			\]
			is lower-semicontinuous.
			From nonnegativity of $u$ it follows that for every $x\in X\setmin K$ and $r \in \R$,
			\begin{equation}
				\label{e:selector1}
				K' \Let \bigl\{a\in A(x) \big| T'u(x, a) \le r\bigr\}\subset\bigl\{a\in A(x)\big| c(x, a) \le r\bigr\},
			\end{equation}
			and the set $\bigl\{a\in A(x)\big| c(x, a) \le r\bigr\}$ is compact by inf-compactness of $c$. Since by definition $T u(x) = \inf_{A(x)} T'u(x, a)$, by Lemma \ref{l:basicselector} it would follow that a selector $f$ exists such that $T u(x) = T'u(x, f(x))\;\;\fa x\in X\setmin K$ once we verify the hypotheses of this Lemma. For this we only have to verify that $T'u$ is l.s.c.\ (which implies it is measurable) and inf-compact on $\mathbb K$. We have seen above that $T'u$ is a l.s.c.\ function on $\mathbb K$. Therefore, for each $x\in X\setmin K$ the map $T'u(x, \cdot)$ is also l.s.c.\ on $A(x)$. Thus, by definition of lower semicontinuity, the set $K'$ in~\eqref{e:selector1} is closed for every $x\in X\setmin K$ and $r\in\R$. Since a closed subset of a compact set is compact, it follows that $K'$ is compact, which in turn shows inf-compactness of $T'u$ on $\mathbb K$ and proves the assertion.
		\end{proof}

		The following lemma shows how functions $u\in\Lp 0(X\setmin K)^+$ satisfying $u\ge T u$ relate to the optimal value function.

		\begin{lemma}
			\label{l:Tineq}
			Suppose that Assumption {\rm \ref{a:key}} holds. If $u\in \Lp 0(X\setmin K)^+$ is such that $u\ge T u$, then $u\ge V^\star$.
		\end{lemma}
		\begin{proof}
			Suppose $u\in\Lp 0(X\setmin K)^+$ satisfies $u \ge T u$, and let $f$ be a selector (whose existence is guaranteed by Lemma \ref{l:selector}) that attains the infimum in~\eqref{e:Tdef}. Fix $x\in X\setmin K$. We have
			\[
				u(x) \ge T u(x) = c(x, f) + \alpha \int_X Q(\mrm d\xi_1|x, f)\:\indic{X\setmin K}(\xi_1)u(\xi_1).
			\]
			The operator $T$ in~\eqref{e:Tdef} is monotone, for if $u, u'\in \Lp 0(X\setmin K)^+$ are two functions with $u \le u'$, then clearly $T u \le T u'$ due to nonnegativity of $c$. Therefore, iterating the above inequality for a second time we obtain
			\begin{equation*}
			\begin{aligned}
				u(x) & \ge c(x, f) + \alpha \int_X Q(\mrm d\xi_1|x, f)\:\indic{X\setmin K}(\xi_1) c(\xi_1, f)\\
				& \qquad + \alpha^2\int_X Q(\mrm d\xi_1|x, f)\:\indic{X\setmin K}(\xi_1)\int_X Q(\xi_2|\xi_1,f)\: \indic{X\setmin K}(\xi_2)u(\xi_2).
			\end{aligned}
			\end{equation*}
			After $n$ such iterations we arrive at
			\[
				u(x) \ge \mathsf E^{f^\infty}_x\!\left[\sum_{i=0}^{(n-1)\mn(\tau-1)} \alpha^i c(x_i, f)\right] + \mathsf E^{f^\infty}_x\bigl[\alpha^{n} u(x_{n})\indic{\{n < \tau\}}\bigr].
			\]
			Since $u\ge 0$, letting $n\ra\infty$ we get
			\[
				u(x) \ge V(f, x) \ge V^\star(x).
			\]
			Since $x\in X\setmin K$ is arbitrary, the assertion follows.
		\end{proof}

		The next lemma deals with convergence of the value iterations to the optimal value function.

		\begin{lemma}
			\label{l:VIconv}
			Suppose that Assumption {\rm \ref{a:key}} holds. Then $v_n\uparrow V^\star$ on $X\setmin K$, and the function $V^\star$ satisfies the $\alpha$-DCOE \eqref{e:alphadcoe}.
		\end{lemma}
		\begin{proof}
			Note that since $v_n(x) = \inf_{\pi\in\Pi} \mathsf E_x^\pi\left[\sum_{i=0}^{(n-1)\mn(\tau-1)} \alpha^i c(x_i, a_i)\right]$ for $x\in X\setmin K$, it follows that 
			\[
				v_n(x) \le \mathsf E_x^\pi\!\left[\sum_{i=0}^{(n-1)\mn(\tau-1)} \alpha^i c(x_i, a_i)\right] \le \mathsf E_x^\pi\!\left[\sum_{i=0}^{\tau-1} \alpha^i c(x_i, a_i)\right],
			\]
			and therefore, taking the infimum over all policies $\pi\in\Pi$ on the right hand side, we get
			\begin{equation}
				\label{e:vnbound}
				v_n(x) \le V^\star(x)\qquad \fa x\in X\setmin K.
			\end{equation}
			Since the cost-per-stage function is nonnegative, $T$ is a monotone operator. Therefore, since $v_0 \Let 0$ and $v_n = T v_{n-1}$ for $n\in\N$, it follows that the $\alpha$-VI functions form a nondecreasing sequence in $\Lp 0(X\setmin K)^+$, which implies that $v_n\uparrow v^\star$ for some function $v^\star \in \Lp 0(X\setmin K)^+$. For $n\in\N$ we define
			\begin{align*}
				\mathbb K\ni (x, a)& \mapsto T'v_n(x, a) \Let c(x, a) + \alpha \int_{X} Q(\mrm dy|x, a)\indic{X\setmin K}(y)v_n(y)\in\R,\\
				\mathbb K\ni (x, a)& \mapsto T'v^\star(x, a) \Let c(x, a) + \alpha \int_{X} Q(\mrm dy|x, a)\indic{X\setmin K}(y) v^\star(y)\in\R.
			\end{align*}
			The monotone convergence theorem guarantees that $T'v_n\uparrow T'v^\star$ pointwise on $\mathbb K$. As in the proof of Lemma~\ref{l:selector} one can establish inf-compactness and lower semicontinuity of $T'v_n$, and $T'v^\star$ on $\mathbb K$. From Lemma~\ref{l:keyconvergence} it now follows that for every $x\in X\setmin K$ we have
			\begin{align*}
				v^\star(x) & = \lim_{n\ra\infty}v_n(x) = \lim_{n\ra\infty} T v_{n-1}(x)\\
				& = \lim_{n\ra\infty} \min_{A(x)}T'v_{n-1}(x, a) = \min_{A(x)} T'v^\star(x, a)\\
				& = T v^\star(x).
			\end{align*}
			This shows that $v^\star$ satisfies the $\alpha$-DCOE, $v^\star = T v^\star$.

			It remains to show that $v^\star = V^\star$. But by Lemma \ref{l:Tineq}, $v^\star = T v^\star$ implies that $v^\star \ge V^\star$ and the reverse inequality follows from~\eqref{e:vnbound} by taking limits as $v^\star = \lim_{n\ra\infty} v_n \le V^\star$.
		\end{proof}

		\begin{lemma}
			\label{l:adcoestat}
			For every deterministic stationary policy $f^\infty$ we have
			\begin{equation}
				\label{e:adcoepolicy}
				V(f^\infty, x) = c(x, f) + \alpha\int_X Q(\mrm dy|x, f)\:\indic{X\setmin K}(y)V(f^\infty, y)\qquad \fa x\in X\setmin K.
			\end{equation}
		\end{lemma}
		\begin{proof}
			Fix a deterministic stationary policy $f^\infty$ and $x\in X\setmin K$. The $\alpha$-discounted cost $V(f^\infty, x)$ corresponding to this policy satisfies, in view of the definition of $\tau$ and the fact that $x\in X\setmin K$,
			\begin{align}
				V(f^\infty, x) & \Let \mathsf E_x^{f^\infty}\!\left[\sum_{i=0}^{\tau-1} \alpha^i c(x_i, f)\right] = \mathsf E_x^{f^\infty}\!\left[ c(x, f) + \sum_{i=1}^{\tau-1} \alpha^i c(x_i, f)\right]\nonumber\\
				& = c(x, f) + \alpha\mathsf E_x^{f^\infty}\!\left[\sum_{i=1}^{\tau-1} \alpha^{i-1} c(x_i, f)\right].\label{e:spolicy1}
			\end{align}
			But then by the Markov property,
			\begin{align*}
				\mathsf E_x^{f^\infty}\!\left[\sum_{i=1}^{\tau-1} \alpha^{i-1} c(x_i, f)\right] & = \mathsf E^{f^\infty}\!\left[\mathsf E^{f^\infty}\!\left[\sum_{i=1}^{\tau-1} \alpha^{i-1} c(x_i, f)\left.\left.\vphantom{\sum_i^\tau}\right|x_{1\mn(\tau-1)}\right]\right|x_0 = x\right]\\
				& = \int_X \indic{X\setmin K}(y) Q(\mrm dy|x, f)\; \mathsf E^{f^\infty}\!\left[\sum_{i=1}^{\tau-1} \alpha^{i-1} c(x_i, f)\left.\vphantom{\sum_i^\tau}\right|x_1 = y\right]\\
				& = \int_X \indic{X\setmin K}(y) Q(\mrm dy|x, f)\; V(f^\infty, y).
			\end{align*}
			This substituted back in~\eqref{e:spolicy1} gives~\eqref{e:adcoepolicy}.
		\end{proof}

		\begin{proof}[Proof of Theorem {\rm \ref{t:EDC}}]
			(i) That $V^\star$ is a solution of the $\alpha$-DCOE follows from Lemma \ref{l:VIconv}, and that $V^\star$ is the minimal solution follows from Lemma~\ref{l:Tineq}, since $u = T u$ implies $u \ge V^\star$.

			(ii) Lemma~\ref{l:selector} guarantees the existence of a selector $f_\star\in\mathbb F$ such that~\eqref{e:alphado} holds. Fix $n\in\N$ and $x\in X\setmin K$. As in the proof of Lemma~\ref{l:Tineq}, iterating equation~\eqref{e:alphado} $n$-times we arrive at
			\begin{align*}
				V^\star(x) & = \mathsf E_x^{f_\star^\infty}\!\left[\sum_{i=0}^{(n-1)\mn(\tau-1)} \alpha^i c(x_i, f_\star)\right] + \mathsf E^{f_\star^\infty}_x\bigl[\alpha^{n} V^\star(x_{n})\indic{\{n < \tau\}}\bigr] \ge \mathsf E_x^{f_\star^\infty}\!\left[\sum_{i=0}^{(n-1)\mn(\tau-1)} \alpha^i c(x_i, f_\star)\right].
			\end{align*}
			By the monotone convergence theorem we have
			\[
				V^\star(x) \ge \lim_{n\ra\infty} \mathsf E_x^{f_\star^\infty}\!\left[\sum_{i=0}^{(n-1)\mn(\tau-1)} \alpha^i c(x_i, f_\star)\right] = \mathsf E_x^{f_\star^\infty}\!\left[\sum_{i=0}^{\tau-1} \alpha^i c(x_i, f_\star)\right],
			\]
			which shows that $V^\star(x) \ge V(f_\star^\infty, x)$, and since $x\in X\setmin K$ is arbitrary, it follows that $V^\star(\cdot) \ge V(f_\star^\infty, \cdot)$. The reverse inequality follows from the definition of $V^\star$ in~\eqref{e:problem}. We conclude that $V^\star(\cdot) = V(f_\star^\infty, \cdot)$, and that $f_\star^\infty$ is an optimal policy.

			For the converse, if $f_\star^\infty$ is an optimal deterministic stationary policy, then by Lemma \ref{l:adcoestat}, equation~\eqref{e:adcoepolicy} becomes
			\[
				V^\star(x) = V(f_\star^\infty, x) = c(x, f_\star) + \alpha \int_X Q(\mrm dy|x, f_\star)\:\indic{X\setmin K}(y)V(f_\star^\infty, y)
			\]
			for $x\in X\setmin K$, which is identical to~\eqref{e:alphado}.
		\end{proof}

	\section{A Contraction Mapping Approach}
	\label{s:contr}
		For the purposes of this section we let $\Lp 0(X\setmin K)$ denote the real vector space of real-valued measurable functions on $X$, and $\Lp 0(X\setmin K)^+$ be the convex cone of nonnegative elements of $\Lp 0(X\setmin K)$. (Note that according to paragraph~\ref{pgr:convention} we let the elements of $\Lp 0(X\setmin K)^+$ take the value $+\infty$.) Given a measurable \emph{weight function} $w:X\setmin K\lra[1, \infty[$ in $\Lp 0(X\setmin K)^+$, we define the weighted norm $\norm{u}_w \Let \sup_{x\in X} \abs{u(x)}/w(x)$. It is well-known that $\bigl(\Lp 0(X\setmin K), \norm{\cdot}_w\bigr)$ is a Banach space.

		\begin{assumption}
			\label{a:further}
			In addition to Assumption~\ref{a:key}, we require that there exist $\ol c > 0$, $\beta\in[1, 1/\alpha[$, and a measurable weight function $w:X\setmin K\lra[1, \infty[$ such that for every $x\in X\setmin K$
			\begin{enumerate}[label=(\roman*), leftmargin=*, align=right, widest=iii]
				\item $\displaystyle{\sup_{A(x)} c(x, a) \le \ol c w(x)}$;
				\item $\displaystyle{\sup_{A(x)} \int_{X\setmin K} Q(\mrm dy|x, a) w(y) \le \beta w(x)}$.\AssumptionEnd
			\end{enumerate}
		\end{assumption}

		\begin{remark}
			If $c$ is bounded, the weight function $w$ may be taken to be $\indic{X\setmin K}$. Also, if $x$ and $x^+$ are the current and the next states of the Markov control process, respectively, then Assumption~\ref{a:further}(ii) implies that 
			\[
				\sup_{A(x)}\mathsf E\bigl[w(x^+)\indic{\{x^+\in X\setmin K\}}\big|(x, a)\bigr] \le \beta w(x)\qquad \fa x\in X\setmin K.
			\]
			We observe that this bears a resemblance with classical Lyapunov-like stability criteria, more specifically, the Foster-Lyapunov conditions~\cite[Chapter~8]{ref:meynCTCN}, \cite{ref:foss04}. However, the condition in Assumption~\ref{a:further}(ii) is uniform over the set of actions $A(x)$ pointwise in $x$. It connects the growth of the cost-per-stage function $c$ with a contraction induced by the discount factor $\alpha$.\RemarkEnd
		\end{remark}

		Recall that a mapping $f:Y\lra Y$ on a nonempty complete metric space $(Y, \rho)$ is a \emph{contraction} if there exists a constant $\gamma\in[0, 1[$ such that $\rho(f(x_1), f(x_2)) \le \gamma\rho(x_1, x_2)$ for all $x_1, x_2\in Y$. The constant $\gamma$ is said to the the \emph{modulus} of the map $f$. A contraction has a unique fixed point $x^\star\in Y$ satisfying $f(x^\star) = x^\star$.

		\begin{proposition}[{\cite[Proposition~7.2.9]{ref:hernandez-lerma2}}]
			\label{p:contr}
			Let $T$ be a monotone map from the Banach space $\bigl(\Lp 0(X\setmin K), \norm{\cdot}_w\bigr)$ into itself. If there exists a $\gamma\in[0, 1[$ such that
			\begin{equation}
				\label{e:contr}
				T(u+rw) \le T(u) + \gamma rw\qquad\text{whenever}\quad u\in\bigl(\Lp 0(X\setmin K), \norm{\cdot}_w\bigr),\quad r\in\R,
			\end{equation}
			then $T$ is a contraction with modulus $\gamma$.
		\end{proposition}

		We have the following lemma.
		\begin{lemma}
			\label{l:Tcontr}
			Under Assumption {\rm \ref{a:further}}, the map $T$ in~\eqref{e:Tdef} is a contraction on $\bigl(\Lp 0(X\setmin K)^+, \norm{\cdot}_w\bigr)$ with modulus $\gamma = \alpha\beta < 1$.
		\end{lemma}
		\begin{proof}
			Fix $u\in\Lp 0(X\setmin K)^+$ with $\norm{u}_w < \infty$. As in the proof of Lemma~\ref{l:selector}, the mapping
			\[
				\mathbb K\ni(x, a)\mapsto T'u(x, a) = c(x, a) + \alpha\int_{X\setmin K} Q(\mrm dy|x, a) u(y)\in\posR
			\]
			is well-defined and l.s.c.\ in $a\in A(x)$ for all $x\in X\setmin K$. By the same Lemma we also know that $T$ maps $\Lp 0(X\setmin K)^+$ into $\Lp 0(X\setmin K)^+$. For every $(x, a)\in\mathbb K$, by Assumption~\ref{a:further},
			\begin{align*}
				\abs{T'(x, a)} & \le c(x, a) + \alpha\int_{X\setmin K}Q(\mrm dy|x, a) \frac{u(y)}{w(y)} w(y) \le \ol cw(x) + \alpha \norm{u}_w\int_{X\setmin K} Q(\mrm dy|x, a) w(y)\\
				& \le \bigl(\ol c + \alpha\beta\norm{u}_w\bigr) w(x),
			\end{align*}
			which shows that $\norm{T'u}_w \le \ol c + \alpha\beta\norm{u}_w$. Therefore, $T$ maps $\bigl(\Lp 0(X\setmin K)^+, \norm{\cdot}_w\bigr)$ into itself. Since $c\ge 0$, it is clear that $T$ is a monotone map on $\bigl(\Lp 0(X\setmin K)^+, \norm{\cdot}_w\bigr)$. By Assumption~\ref{a:further}(ii), for $r\in\R$ and $x\in X\setmin K$ we have
			\begin{align*}
				T(u+rw)(x) & = \min_{A(x)}\left(c(x, a) + \alpha\int_{X\setmin K} Q(\mrm dy|x, a)\bigl(u(y) + rw(y)\bigr)\right)\\
					& \le \min_{A(x)}\left(c(x, a) + \alpha\int_{X\setmin K} Q(\mrm dy|x, a)u(y)\right) + r\alpha\beta w(x)\\
					%& \qquad\qquad + \alpha r\left.\indic{X\setmin K}(x)\int_{X\setmin K} Q(\mrm dy|x, a) w(y)\right)\\
					& \le Tu(x) + r\alpha\beta w(x).
			\end{align*}
			This shows that~\eqref{e:contr} holds with $\gamma = \alpha\beta$, and Proposition~\ref{p:contr} implies that $T$ is a contraction on $\bigl(\Lp 0(X\setmin K)^+, \norm{\cdot}_w\bigr)$.
		\end{proof}

		The following proposition establishes bounds for the distance between the optimal value function $V^\star$ and the $\alpha$-VI functions $(v_n)_{n\in\Nz}$ by employing the contraction mapping $T$ of Lemma~\ref{l:Tcontr}.

		\begin{proposition}
			\label{p:Tfp}
			Suppose that Assumption {\rm \ref{a:further}} holds, and let $\gamma \Let \alpha\beta$. Then:
			\begin{enumerate}[label=\emph{(\roman*)}, align=right, leftmargin=*, widest=iii]
				\item The $\alpha$-discount value function $V^\star$ satisfies $\norm{V^\star}_w \le \ol c/(1-\gamma)$.
				\item The $\alpha$-VI functions $(v_n)_{n\in\Nz}$ satisfy
				\[
					V^\star(x) - v_n(x) \le \ol c w(x)\left(\frac{\gamma^n}{1-\gamma}\right)\qquad \fa x\in X\setmin K, \quad \fa n\in\N.
				\]
				In particular, $\norm{v_n - V^\star}_w \le \ol c\gamma^n/(1-\gamma)\;\;\; \fa n\in\Nz$.
				\item The optimal value function $V^\star$ is the unique function in $\bigl(\Lp 0(X\setmin K)^+, \norm{\cdot}_w\bigr)$ that solves the $\alpha$-DCOE~\eqref{e:alphadcoe}.
			\end{enumerate}
		\end{proposition}
		\begin{proof}
			(i) Let $\pi$ be an arbitrary Markov policy. Trivially we have $\mathsf E^\pi_x\bigl[w(x_0)\bigr] \le w(x)$. Fix $i\in\N$, and a history $h_i\in\sigalg_{i\mn\tau}$. In view of Assumption~\ref{a:further}(ii), on the event $\{i < \tau\}$ we have
			\begin{align*}
				\mathsf E^\pi_x\bigl[w(x_i)\big|h_{i-1}, a_{i-1}\bigr] = \int_{X\setmin K} Q(\mrm dy|x_{i-1}, a_{i-1}) w(y) \le \beta w(x_{i-1})\quad \fa a_i\in A(x_i),
			\end{align*}
			which shows that $\mathsf E^\pi_x\bigl[w(x_i)\indic{\{i < \tau\}}\bigr] \le \beta \mathsf E^\pi_x\bigl[w(x_{i-1})\indic{\{i < \tau\}}\bigr]$. Iterating this inequality we arrive at $\mathsf E^\pi_x\bigl[w(x_i)\indic{\{i < \tau\}}\bigr] \le \beta^i w(x)$.  Also, by Assumption~\ref{a:further}(i) we have $c(x_i, a_i) \le \ol c w(x_i)$ for all $i\in\Nz$ such that $i < \tau$, which in conjunction with the above inequality gives
			\begin{equation}
			\label{e:Tfp1}
				\mathsf E^\pi_x\bigl[c(x_i, a_i)\indic{\{i < \tau\}}\bigr] \le \ol c\beta^i w(x).
			\end{equation}
			By the monotone convergence theorem and~\eqref{e:Tfp1} we have
			\begin{equation}
			\label{e:Tfp2}
			\begin{aligned}
				V(\pi, x) & = \mathsf E^\pi_x\!\left[\sum_{i=0}^\infty \alpha^i c(x_i, a_i) \indic{\{i < \tau\}}\right] \le \sum_{i=0}^\infty \alpha^i \mathsf E^\pi_x\bigl[c(x_i, a_i)\indic{\{i < \tau\}}\bigr]\\
				& \le \ol c\sum_{i=0}^\infty (\alpha\beta)^i w(x) \le w(x)\cdot\frac{\ol c}{1-\gamma}.
			\end{aligned}
			\end{equation}
			It follows immediately that $\norm{V^\star}_w = \norm{\inf_{\Pi} V(\pi, x)}_w \le \ol c/(1-\gamma)$.

			(ii) By definition, the $\alpha$-VI functions $(v_n)_{n\in\Nz}$ satisfy $v_n = T v_{n-1} = T^n v_0$, with $v_0 \Let 0$. Since $T$ is a contraction on $\bigl(\Lp 0(X\setmin K)^+, \norm{\cdot}_w\bigr)$ by Lemma \ref{l:Tcontr}, it follows that $T$ has a unique fixed point, which, by definition is $V^\star$, since $\norm{V^\star}_w < \infty$ by (i). A standard property of contraction maps implies that
			\[
				\norm{T^n v_0 - V^\star}_w \le \gamma^n\norm{v_0 - V^\star}_w \qquad\fa u\in\Lp 0(X\setmin K)^+, \norm{u}_w < \infty,\quad \fa n\in\Nz.
			\]
			With the bound on $\norm{V^\star}_w$ obtained in (i), we get $\norm{v_n - V^\star}_w \le \ol c\cdot\gamma^n/(1-\gamma)$. Since $T$ is also a contraction on $\bigl(\Lp 0(X\setmin K)^+, \norm{\cdot}_w\bigr)$, $v_n|_{K} = 0$, and $v_n\uparrow V^\star$, the last inequality yields $V^\star(x) - v_n(x) \le \ol cw(x)\gamma^n/(1-\gamma)$ for every $x\in X\setmin K$.

			(iii) Of course $V^\star$ solves the $\alpha$-DCOE~\eqref{e:alphadcoe}. Uniqueness follows from the facts that the operator $T$ in~\eqref{e:Tdef} is a contraction by Lemma \ref{l:Tcontr}, and that the fixed point of a contraction mapping in a Banach space (or more generally, in a complete metric space) is unique.
		\end{proof}

		Note that the conditions in Assumption~\ref{a:further} are automatic if $c$ is bounded. This gives the following straightforward result.
		\begin{corollary}
			\label{c:Tfp}
			Suppose that Assumption {\rm \ref{a:key}} holds, and $\wt c \Let \sup_{\mathbb K}c(x, a) < \infty$. Then:
			\begin{enumerate}[label=\emph{(\roman*)}, align=right, leftmargin=*, widest=iii]
				\item The $\alpha$-discount value function $V^\star$ satisfies $\norm{V^\star} \le \wt c/(1-\alpha)$.
				\item The $\alpha$-VI functions $(v_n)_{n\in\Nz}$ satisfy
				\[
					V^\star(x) - v_n(x) \le \wt c \left(\frac{\alpha^n}{1-\alpha}\right)\qquad \fa x\in X\setmin K,\quad \fa n\in\N.
				\]
				In particular, $\norm{v_n - V^\star} \le \wt c\:\alpha^n/(1-\alpha)\;\;\; \fa n\in\Nz$.
				\item The optimal value function $V^\star$ is the unique function in $\bigl(\Lp 0(X\setmin K)^+, \norm{\cdot}_w\bigr)$ that solves the $\alpha$-DCOE~\eqref{e:alphadcoe}.
			\end{enumerate}
		\end{corollary}

	\section{Asymptotic Discount Optimality of the $\alpha$-VI Policy}
	\label{s:ado}
		We have seen that the $\alpha$-value iteration functions $(v_n)_{n\in\Nz}$ defined in~\eqref{e:VI} converge to $V^\star$ by Lemma~\ref{l:VIconv}. In this section we address the question whether the $\alpha$-VI policies converge in some sense to a policy $f_\star^\infty$ as $n\ra\infty$.

		\begin{defn}
			\label{d:avip}
			Let $(v_n)_{n\in\Nz}$ be the sequence of $\alpha$-VI functions in~\eqref{e:VI}, and let $\wh\pi = (\wh f_n)_{n\in\Nz}\in\Pi_{DM}$ be a deterministic Markov policy such that $\wh f_0\in\mathbb F$ is arbitrary, and for $n\in\N$,
			\[
				v_n(x) = c\bigl(x, \wh f_n\bigr) + \alpha \int_{X\setmin K}Q\bigl(\mrm dy\big|x, \wh f_n\bigr) v_{n-1}(y)\qquad \fa x\in X\setmin K.
			\]
			Then $\wh\pi$ is called an \emph{$\alpha$-VI policy}.\DefEnd
		\end{defn}

		Under Assumption~\ref{a:key} we get the following basic existential result.

		\begin{proposition}
			\label{p:ado}
			Suppose that Assumption {\rm \ref{a:key}} holds, the action space $A$ is locally compact, and let $\wh\pi = \bigl(\wh f_n\bigr)_{n\in\Nz}\in\Pi_{DM}$ be an $\alpha$-VI policy as defined in Definition~\ref{d:avip}. Then there exists a selector $\wh f\in\mathbb F$ such that for every $x\in X\setmin K$, $\wh f(x)\in A(x)$ is an accumulation point of $\bigl(\wh f_n(x)\bigr)_{n\in\Nz}$, and the corresponding deterministic stationary policy $\wh f^\infty\in\Pi_{DS}$ is $\alpha$-discount optimal.
		\end{proposition}

		The proof is based on the following immediate adaptation of \cite[Lemma 4.6.6]{ref:hernandez-lerma1}.
		\begin{lemma}
			\label{l:ado}
			Let $u$ and $u_n$, $n\in\N$, be l.s.c.\ functions, bounded below, and inf-compact on $\mathbb K$. For every $n\in\N$ let $u_n^\star(x) \Let \min_{A(x)} u_n(x, a)$ and $u^\star(x) \Let \min_{A(x)} u(x, a)$, let $\wh f_n\in\mathbb F$ be a selector such that $u_n^\star(x) = u_n\bigl(x, \wh f_n(x)\bigr)$ for all $x\in X\setmin K$. If $A$ is locally compact and $u_n\uparrow u$, then there exists a selector $\wh f\in\mathbb F$ such that $\wh f(x)\in A(x)$ is an accumulation point of the sequence $\bigl(\wh f_n(x)\bigr)_{n\in\N}$ for every $x\in X\setmin K$, and $u^\star(x) = u\bigl(x, \wh f(x)\bigr)$.
		\end{lemma}

		\begin{proof}[Proof of Proposition {\rm \ref{p:ado}}]
			For $(x, a)\in\mathbb K$ we define $u(x, a) \Let c(x, a) + \alpha \int_{X\setmin K} Q(\mrm dy|x, a) V^\star(y)$, and
			\begin{equation}
				\label{e:unconv}
				u_n(x, a) \Let c(x, a) + \alpha \int_{X\setmin K} Q(\mrm dy|x, a) v_{n-1}(y).
			\end{equation}
			Since $c\ge 0$, the functions $u_n$ and $u$ are nonnegative. Since $v_n\uparrow V^\star$ by Lemma~\ref{l:VIconv}, the monotone convergence theorem implies that 
			\[
				\int_{X\setmin K}Q(\mrm dy|x, a) v_n(y) \lra \int_{X\setmin K} Q(\mrm dy|x, a) V^\star(y)
			\]
			pointwise on $\mathbb K$. It is clear that $u_n\uparrow u$, and the assertion follows at once from Lemma~\ref{l:ado}.
		\end{proof}

		Under the stronger Assumption~\ref{a:further} we get quantitative estimates of the rate at which the $\alpha$-VI policy defined in Definition~\ref{d:avip} converges to an optimal one.

		\begin{defn}
			The function $D:\mathbb K\lra\posR$ defined by
			\[
				\mathbb K\ni (x, a) \mapsto D(x, a) \Let c(x, a) + \alpha\int_{X\setmin K} Q(\mrm dy|x, a) V^\star(y) - V^\star(x)
			\]
			is called the \emph{$\alpha$-discount discrepancy function}. The $\alpha$-VI policy $\wh\pi = \bigl(\wh f_n\bigr)_{n\in\Nz}$ defined in Definition~\ref{d:avip} is called \emph{pointwise asymptotically discount optimal} if for every $x\in X\setmin K$ we have $\lim_{n\ra\infty} D\bigl(x, \wh f_n\bigr) = 0$.\DefEnd
		\end{defn}

		It is clear that for $x\in X\setmin K$ and a selector $f\in\mathbb F$ (see paragraph~\ref{pgr:convention}), the $\alpha$-discount discrepancy function $D(x, f(x))$ is $0$ if and only if $f^\infty$ is an optimal policy. The function $D$ measures closeness to an optimal selector in a weak sense.

		\begin{proposition}
			Suppose that Assumption {\rm \ref{a:further}} holds, and let $\gamma \Let \alpha\beta$. Then the $\alpha$-VI policy $\wh\pi = \bigl(\wh f_n\bigr)_{n\in\Nz}$ is pointwise asymptotically discount optimal, and for every $x\in X\setmin K$ and $n\in\N$,
			\[
				0 \le D\bigl(x, \wh f_n\bigr) \le 2\ol c\left(\frac{\gamma^{n+1}}{1-\gamma}\right) w(x).
			\]
		\end{proposition}
		\begin{proof}
			The first inequality follows directly from the definition of $V^\star$. To prove the second inequality fix $x\in X\setmin K$. We see that by the definition of the discrepancy function,
			\begin{equation}
			\label{e:adiscop1}
			\begin{aligned}
				D\bigl(x, \wh f_n\bigr) & = c\bigl(x, \wh f_n\bigr) + \alpha \int_{X\setmin K} Q\bigl(\mrm dy\big|x, \wh f_n\bigr) V^\star(y) - V^\star(x)\\
				& = \bigl(v_{n+1}(x) - V^\star(x)\bigr) + \alpha\int_{X\setmin K} Q\bigl(\mrm dy\big|x, \wh f_n\bigr)\bigl(V^\star(y) - v_{n}(y)\bigr).
			\end{aligned}
			\end{equation}
			By Proposition~\ref{p:Tfp}(ii) we have 
			\begin{equation}
			\label{e:adiscop2}
				\abs{v_{n+1}(x) - V^\star(x)} \le \ol cw(x) \frac{\gamma^{n+1}}{1-\gamma},
			\end{equation}
			and in the light of Assumption~\ref{a:further}(ii) we arrive at
			\begin{equation}
			\label{e:adiscop3}
			\begin{aligned}
				\int_{X\setmin K}Q\bigl(\mrm dy\big|x, \wh f_n\bigr)\bigl(V^\star(y) - v_{n}(y)\bigr) & \le \int_{X\setmin K}Q\bigl(\mrm dy\big|x, \wh f_n\bigr)\bigl(V^\star(y) - v_{n}(y)\bigr)\\
				& \le \frac{\gamma^{n}}{1-\gamma}\beta w(x).
			\end{aligned}
			\end{equation}
			The assertion follows immediately after substituting~\eqref{e:adiscop2} and~\eqref{e:adiscop3} in~\eqref{e:adiscop1}.
		\end{proof}

		For bounded costs we have the following straightforward conclusion.

		\begin{corollary}
			Suppose that Assumption {\rm \ref{a:key}} holds, and $\wt c \Let \sup_{\mathbb K}c(x, a) < \infty$. Then the $\alpha$-VI policy $\wh\pi = \bigl(\wh f_n\bigr)_{n\in\Nz}$ is pointwise asymptotically discount optimal, and for every $x\in X\setmin K$ and $n\in\N$,
			\[
				0 \le D\bigl(x, \wh f_n\bigr) \le 2\wt c\:\left(\frac{\alpha^{n+1}}{1-\alpha}\right).
			\]
		\end{corollary}

	\section{Average cost of recovery}
		As mentioned in \secref{s:intro}, a motivation for this work was to come up with a suitable recovery strategy for MPC. Tracing our development of the MPC methodology in \secref{s:intro}, one sees that in the presence of state and/or action constraints, one seeks a deterministic stationary policy $g_\star^\infty$ that is active whenever the state is inside the safe set $K$, and a recovery strategy outside $K$. Let us assume that for a given problem we have determined such a policy, and we have also determined a deterministic stationary policy $f_\star^\infty$ corresponding to the recovery strategy corresponding to a cost-per-stage function defined on $X\setmin K$ for the same problem as described in the preceding sections. One of the natural questions at this stage is whether one can find estimates of the average cost of recovery.%To wit, what is the average cost of the excursions of the state of the process outside $K$ that is incurred when the policies $g_\star^\infty$ and $f_\star^\infty$ are applied inside and outside $K$, respectively?

		To this end let us define two constants:
		\begin{equation}
		\label{e:betadef}
			\beta_1 \Let \inf_{x\in K}\int_{X\setmin K} Q(\mrm dy|x, g_{\star}) V^\star(y),\quad  \beta_2 \Let \sup_{x\in K}\int_{X\setmin K} Q(\mrm dy|x, g_{\star}) V^\star(y),
		\end{equation}
		where $V^\star$ is as defined in \eqref{e:problem}. Let $\wt f^\infty$ be the deterministic stationary policy defined by
		\begin{equation}
		\label{e:overallpolicy}
			\wt f(x) \Let f_\star(x)\indic{X\setmin K}(x) + g_\star(x)\indic{K}(x);
		\end{equation}
		to wit, $\wt f^\infty$ consists of concatenation of $f_\star^\infty$ and $g_\star^\infty$ between exit and entry times to $K$. We have the following result:

		\begin{proposition}
			Let $g_\star^\infty$ be a deterministic stationary policy that is active whenever the state is inside the set $K$, and let $f_\star^\infty$ be a recovery strategy corresponding to the problem \eqref{e:problem}. Let the initial condition $x$ be in $X\setmin K$. We define the average cost of recovery
			\[
				\wt W(x) \Let \lim_{n\to\infty} \frac{1}{n+1} \EE^{\wt f^\infty}_x\Biggl[\sum_{i=0}^{n}\sum_{t=\tau_{2i}}^{\tau_{2i+1}-1} \alpha^{t-\tau_{2i}} c(x_t, a_t)\Biggr],
			\]
			where $\wt f^\infty$ is as defined in \eqref{e:overallpolicy}, $\tau_0 \Let 0$, $\tau_1$ is the first entry time to $K$, $\tau_2$ is the first exit time from $K$ after $\tau_1$, and so on. Suppose that from any initial condition in $X\setmin K$ the first hitting time of $K$ is finite almost surely under $f_\star^\infty$, and from any initial condition in $K$ the first hitting time of $X\setmin K$ is finite almost surely under $g_\star^\infty$. Then we have $\beta_1 \le \wt W(x) \le \beta_2$, where $\beta_1, \beta_2$ are as defined in \eqref{e:betadef}.
		\end{proposition}
		Note that an identical bound holds if the initial condition $x\in K$, with an obvious relabelling of the stopping times $(\tau_i)_{i\in\Nz}$.
		\begin{proof}
			First of all, note that the policy $\wt f^\infty$ is deterministic stationary, and under this policy the controlled process is stationary Markov. Now we have for a fixed $n\in \N$:
			\begin{equation}
			\label{e:costofrecovery}
			\begin{aligned}
				\EE^{\wt f^\infty}_x & \Biggl[\sum_{i=0}^n \sum_{t=\tau_{2i}}^{\tau_{2i+1}-1} \alpha^{t-\tau_{2i}} c(x_t, a_t)\Biggr] = \sum_{i=0}^n \EE^{\wt f^\infty}_x\Biggl[\sum_{t=\tau_{2i}}^{\tau_{2i+1}-1} \alpha^{t-\tau_{2i}} c(x_t, a_t)\Biggr]\\
				& = \sum_{i=0}^n \EE^{\wt f^\infty}_x\Biggl[\EE^{f_\star^\infty}\Biggl[\sum_{t=\tau_{2i}}^{\tau_{2i+1}-1} \alpha^{t-\tau_{2i}} c(x_t, a_t)\Bigg|\sigalg_{\tau_{2i}}\Biggr]\Biggr] = \sum_{i=0}^n \EE^{\wt f^\infty}_x\Bigl[\EE^{f_\star^\infty}\bigl[V^\star(x_{\tau_{2i}})\big|\sigalg_{\tau_{2i}}\bigr]\Bigr],
			\end{aligned}
			\end{equation}
			where the first equality follows from monotone convergence and the last equality from the strong Markov property. Appealing to the strong Markov property once again we see that $\EE^{f_\star^\infty}\bigl[V^\star(x_{\tau_{2i}})\big|\sigalg_{\tau_{2i}}\bigr] = \EE^{f_\star^\infty}\bigl[V^\star(x_{\tau_{2i}})\big|x_{\tau_{2i}}\bigr]$. Finally, from the definition of $\tau_{2i}$ it follows that
			\[
				\EE^{f_\star^\infty}\bigl[V^\star(x_{\tau_{2i}})\big|x_{\tau_{2i}}\bigr] \le \sup_{\xi\in K}\int_X Q(\mrm dy|\xi, f_{\mathrm{in}}) V^\star(y)\indic{X\setmin K}(y) = \beta_2.
			\]
			It is not difficult to arrive at the lower bound $\EE^{f_\star^\infty}\bigl[V^\star(x_{\tau_{2i}})\big|x_{\tau_{2i}}\bigr]\ge \beta_1$ by following the same steps as above. Substituting in~\eqref{e:costofrecovery} and taking limits we arrive at the assertion.
		\end{proof}

	\section{A Rolling Horizon Implementation}
	\label{s:rh}
		The \emph{rolling-horizon} procedure can be briefly described as follows. Fix a horizon $N\in\N$ and set $n = 0$. Then 
		\begin{enumerate}[label=(\alph*), leftmargin=*, align=right]
			\item we determine an optimal control policy, say $\pi^\star_{n:n+N}$, for the $(N+1)$-period cost function starting from time $n$, given the (perfectly observed) initial condition $x_n$; standard arguments lead to a realization of this policy as a sequence of $(N+1)$ selectors $\bigl\{\wh f_{n, n+N-j}\big|j=n, n+1, \ldots, n+N\bigr\}$;% and we define $\wh f_N \Let f_{n, N}$;
			\item we increase $n$ to $n+1$, and go back to step (a).
		\end{enumerate}
		Accordingly, the $n$-th step of this procedure consists of minimizing the stopped $(N+1)$-period cost function starting at time $n$, namely, the objective is to find a control policy that attains
		\begin{equation}
		\label{e:rhcf}
		\begin{aligned}
			\inf_{\pi\in\Pi} V_{n, n+N}(\pi, x) \Let \;\inf_{\pi\in\Pi}\mathsf E^\pi\!\left[\left.\sum_{i=n\mn(\tau-1)}^{(n+N)\mn(\tau-1)} \alpha^{i-n\mn(\tau-1)} c(x_i, a_i)\right|x_{n\mn(\tau-1)} = x\right]
		\end{aligned}
		\end{equation}
		for $x\in X\setmin K$. By stationarity and Markovian nature of the control model, it is enough to consider the control problem of minimizing the cost for $n = 0$, i.e., the problem of minimizing $V_{0, N}(\pi, x)$ over $\pi\in\Pi$. The corresponding policy $\pi$ is given by the policy that minimizes the $(N+1)$-stage $\alpha$-VI function $v_{N+1}$ in~\eqref{e:VI}. This particular policy is realized as a sequence of $(N+1)$ selectors $\bigl(\wh f_N, \ldots, \wh f_0\bigr)$. Thus, in the light of the above discussion, the rolling-horizon procedure yields the stationary suboptimal control policy $\wh\pi \Let \wh f_N^\infty$ for the original problem~\eqref{e:problem}.

		Let $V\bigl(\wh f_N^\infty, x\bigr)$ be the value function corresponding to the deterministic stationary policy $\wh f_N^\infty \Let \bigl(\wh f_N, \wh f_N, \ldots\bigr)$, $x\in X\setmin K$. Observe that $\norm{V\bigl(\wh f_N^\infty, x\bigr)}_w < \infty$, which follows from the more general estimate in~\eqref{e:Tfp2}. Our objective in this section is to give quantitative estimates of the extent of sub-optimality of the rolling-horizon policy $\wh\pi$, compared to the optimal policy $\pi^\star$ that attains the infimum in~\eqref{e:problem}. We shall follow the notations of~\secref{s:contr} above.

		\begin{theorem}
			\label{t:rh}
			Suppose that Assumption {\rm \ref{a:further}} holds, and let $\gamma \Let \alpha\beta$. For every $N\in\Nz$ and $x\in X\setmin K$ we have
			\begin{equation}
			\label{e:keyrh}
				0 \le V\bigl(\wh f_N^\infty, x\bigr) - v_{N+1}(x) \le \ol c w(x) \left(\frac{\gamma^{N+1}}{1-\gamma}\right),
			\end{equation}
			where $v_{N+1}$ is the $(N+1)$-th $\alpha$-VI function defined in~\eqref{e:VI}. In particular, 
			\begin{equation}
			\label{e:sloppyrh}
				V\bigl(\wh f_N^\infty, x\bigr) - V^\star(x) \le \ol c w(x) \left(\frac{\gamma^{N+1}}{1-\gamma}\right).
			\end{equation}
		\end{theorem}

		A proof of Theorem \ref{t:rh} is given in the Appendix, if follows the arguments in~\cite[Theorem~1]{ref:aldenRH} for finite state-space Markov decision processes and bounded costs. It is of interest to note that the bound in~\eqref{e:keyrh} is identical to the bound between $V^\star(\cdot)$ and $v_{N+1}(\cdot)$ that appears in Proposition~\ref{p:Tfp}.

		If the cost-per-stage function $c$ is bounded on $\mathbb K$, we have the following immediate corollary:
		\begin{corollary}
			Suppose the Markov control process satisfies Assumption {\rm \ref{a:key}}. Let the cost-per-stage function $c:\mathbb K\lra\posR$ be bounded, with $\wt c \Let \sup_{\mathbb K} c(x, a) < \infty$. Then $V\bigl(\wh f_N^\infty, x\bigr) \ge V^\star(x)$ for every $x\in X\setmin K$, and
			\[
				\sup_{x\in X\setmin K}\left(V\bigl(\wh f_N^\infty, x) - V^\star(x)\right) \le \frac{\wt c\cdot\alpha^{N+1}}{1-\alpha}.
			\]
		\end{corollary}

	\section{Application}
	\label{s:appl}
		In this section we give a numerical example concerning fishery management. The example is motivated by~\cite[Chapter~7]{hastings1989introduction}. The example considers a fishery modeled in discrete-time with the time period representing a fishing season. The state of the controlled Markov chain is the population of the fish species of interest. Fishermen might on the one hand want to harvest all that they can manage in order to increase their short-run profit, but on the other hand this might lead to very low levels of the population.  Our goal is to design a recovery strategy for the case that the population gets over-fished and goes below a critical level.

        For doing so, we consider a simple model, with four possible fish population levels, 1 (almost extinct), 2, 3, and 4 (the target set). We assume that we can accurately measure the population size at the beginning of each season $k$, $X_k$. During a season the following set of actions are available: {Harvest (1), Harvest less (2), Do nothing (3), Import fish (4), Import less (5)}. We also take as given the following transition probabilities between the Markov States, where $T_a(i,j)$ denotes the probability that the population level at the beginning of the next season will be $j$, given that the current population is $i$ and action $a$ is applied during this season.

        \begin{alignat*}{2}
        T_1 & = \begin{bmatrix} 1 & 0 & 0 &0\\
        0.7 & 0.3 & 0 & 0\\
        0.1 & 0.6 & 0.3 & 0\\
        - & - & - & - \end{bmatrix} & \quad &
        T_2 = \begin{bmatrix} 1 & 0 & 0 & 0\\
        0.35 & 0.65 & 0 & 0\\
        0.04 & 0.5 & 0.46 & 0\\
        - & - & - & - \end{bmatrix}\\
        T_3 & = \begin{bmatrix} 0.99 & 0.01 & 0 &0\\
        0.01 & 0.7 & 0.28 & 0.01\\
        0 & 0.03 & 0.65 & 0.32\\
        - & - & - & - \end{bmatrix} & \quad &
        T_4 = \begin{bmatrix} 0.4 & 0.6 & 0 & 0\\
        0 & 0.3 & 0.65 & 0.05\\
        0 & 0 & 0.25 & 0.75\\
        - & - & - & - \end{bmatrix}\\
        T_5 & = \begin{bmatrix} 0.6 & 0.4 & 0 & 0\\
        0 & 0.45 & 0.54 & 0.01\\
        0 & 0 & 0.45 & 0.55\\
        - & - & - & - \end{bmatrix} & &
        \end{alignat*}

        The costs incurred at each state are $c(x_i,\alpha_i) = C(x_i) + A(x_i,\alpha_i)$, where \begin{equation*}C(x_i)=\begin{bmatrix} 300 & 150 & 100 & - \end{bmatrix}^\textrm{T} \end{equation*} represents a cost incurred for being at the current state and \begin{equation*}A(x_i,\alpha_i) = \begin{bmatrix} -20 & -10 & 0 & 150 & 75\\ -40 & -20 & 0 & 150 & 75\\ -80 & -40 & 0 & 150 & 75\\ - & - & - & - & - \end{bmatrix}\end{equation*} the action cost associated with each action and state. We assume a discount factor $\alpha = 0.9$.

        Using this setting, one can compute the policy that attains the $\alpha$-discount value function~\eqref{e:problem}. This turns out to be to import fish when in state $(1)$, to import fewer fish in state $(2)$, and do nothing at state $(3)$. Next, we search for the optimum policy, while using a rolling horizon control scheme, i.e., finding the policy that attains~\eqref{e:rhcf}. We solve the problem for horizon lengths between $1$ and $10$, in order to compare the results with the infinite horizon optimal policy.

        \begin{figure}[h]
          \centering
          \includegraphics[width=0.49\textwidth]{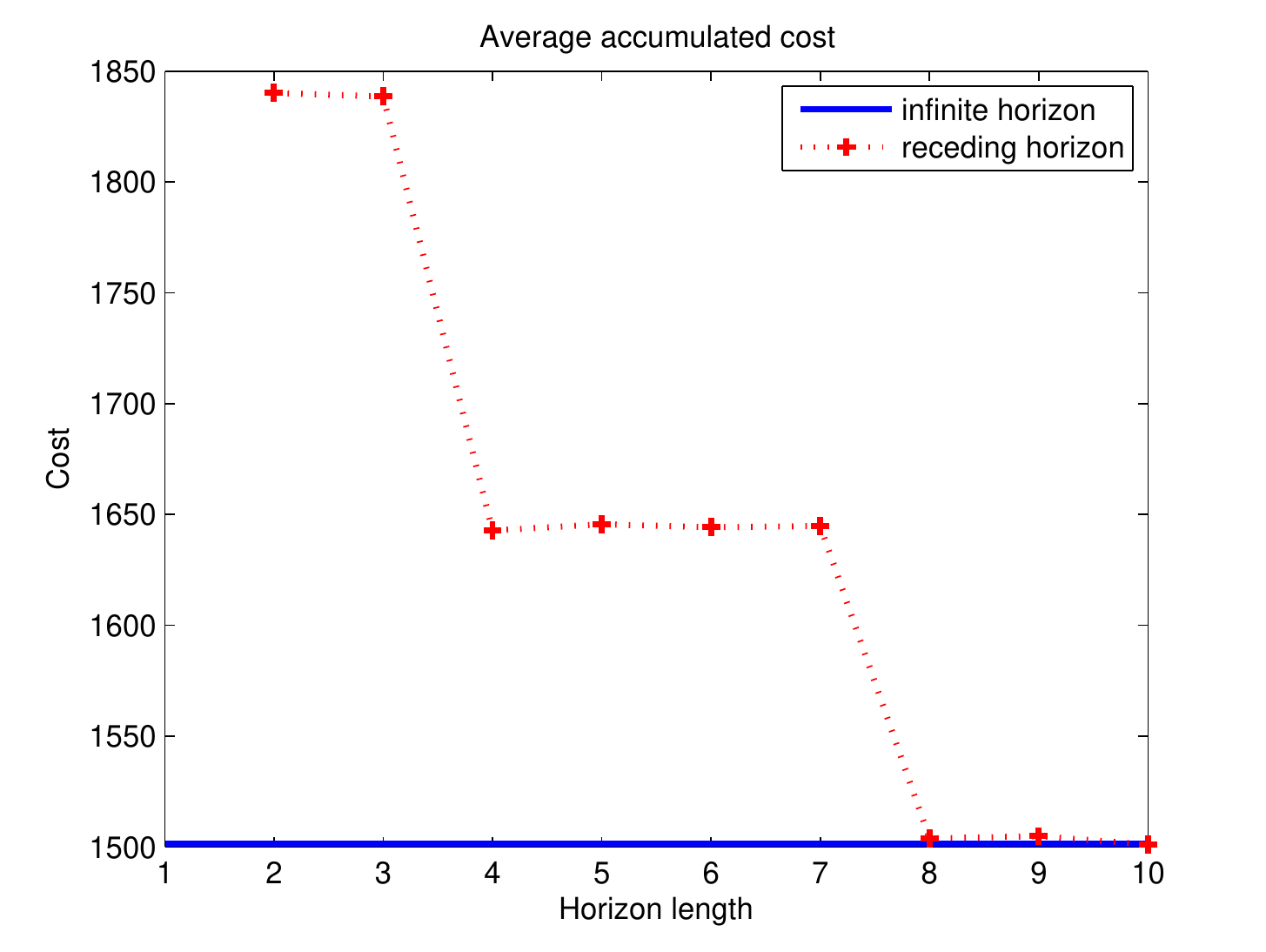}
          \includegraphics[width=0.49\textwidth]{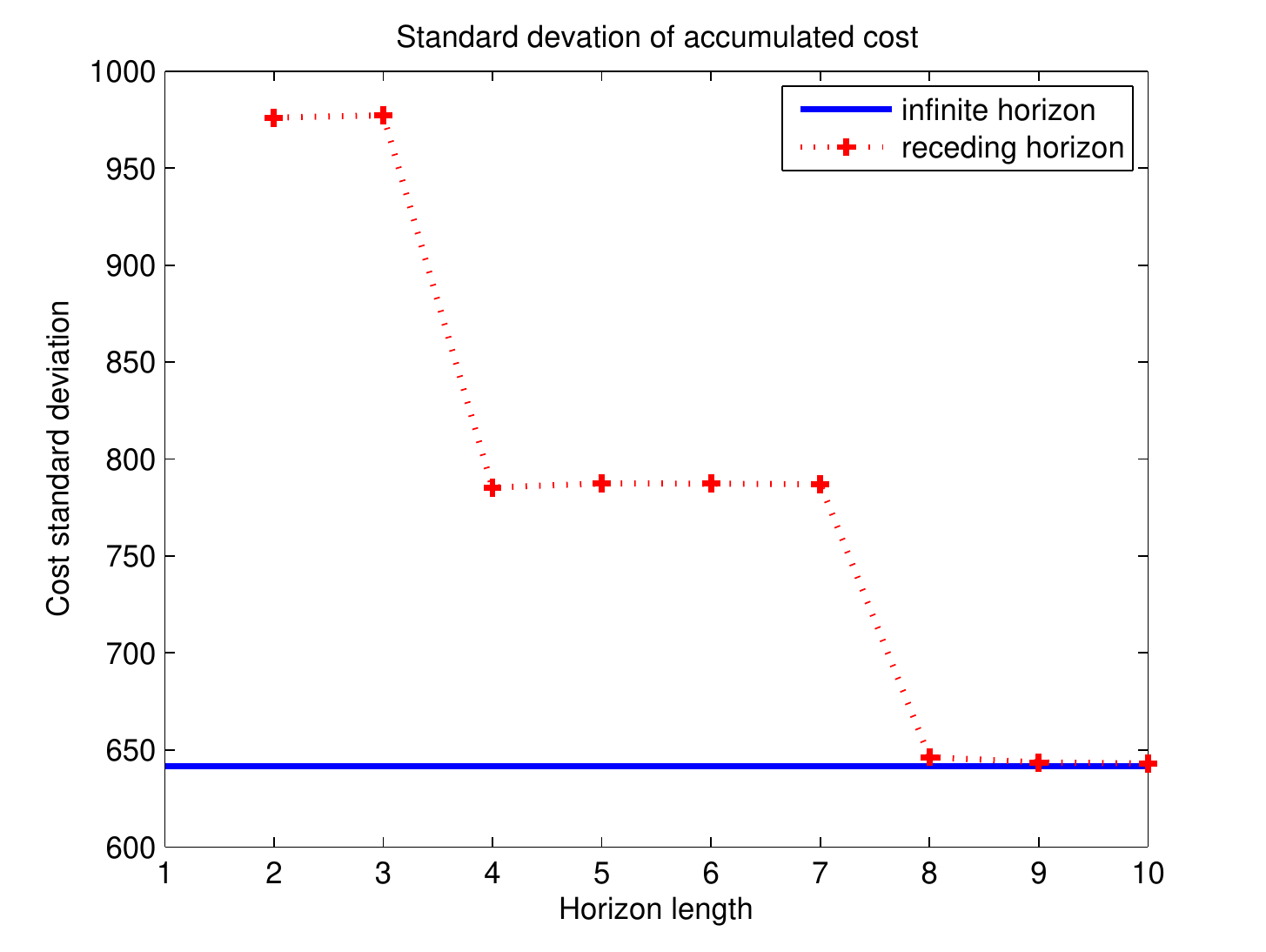}
          \caption{Accumulated cost average and standard-deviation}
          \label{f:cost}
        \end{figure}

        \begin{figure}[h]
          \centering
          \includegraphics[width=0.49\textwidth]{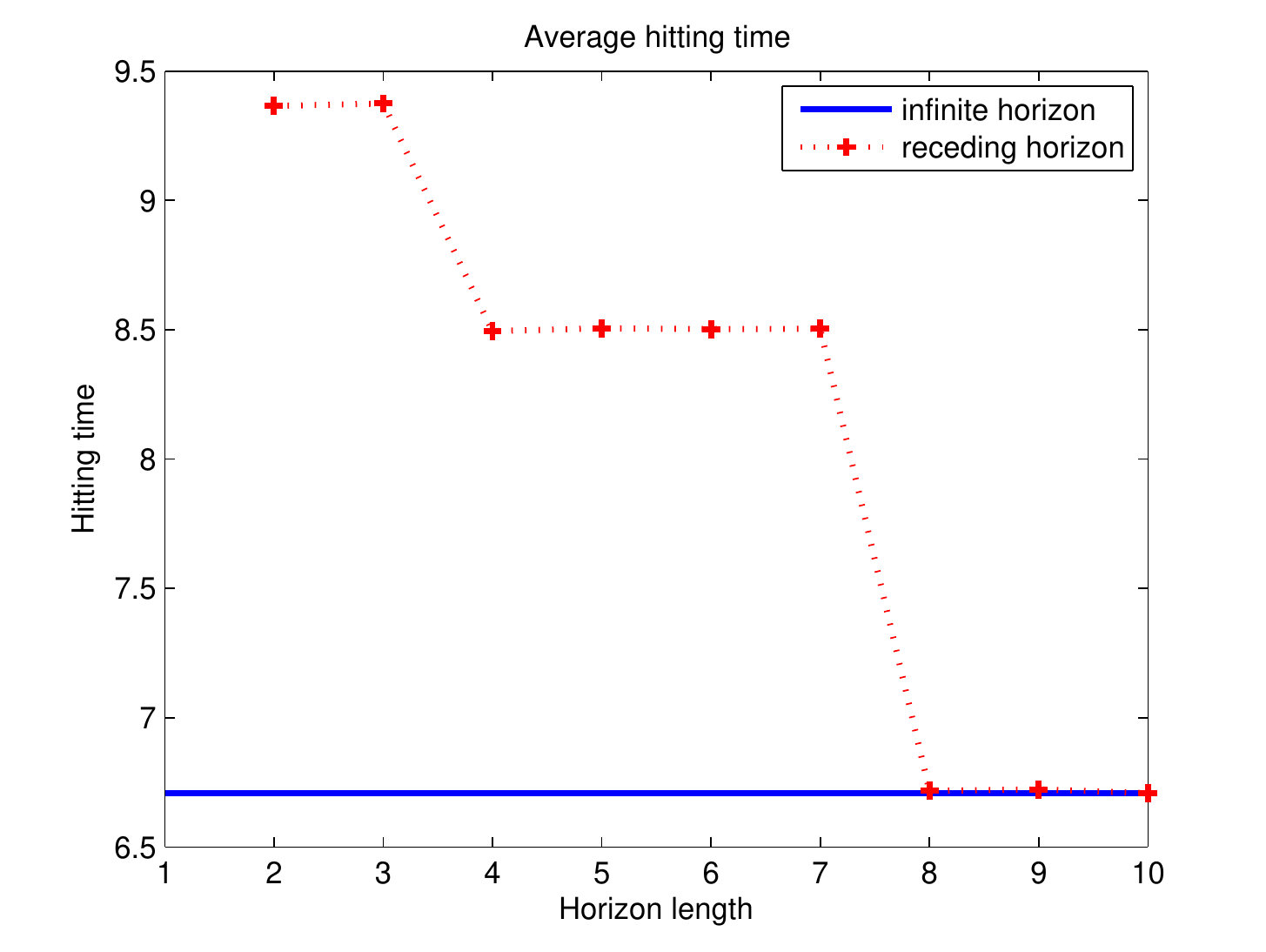}
          \includegraphics[width=0.49\textwidth]{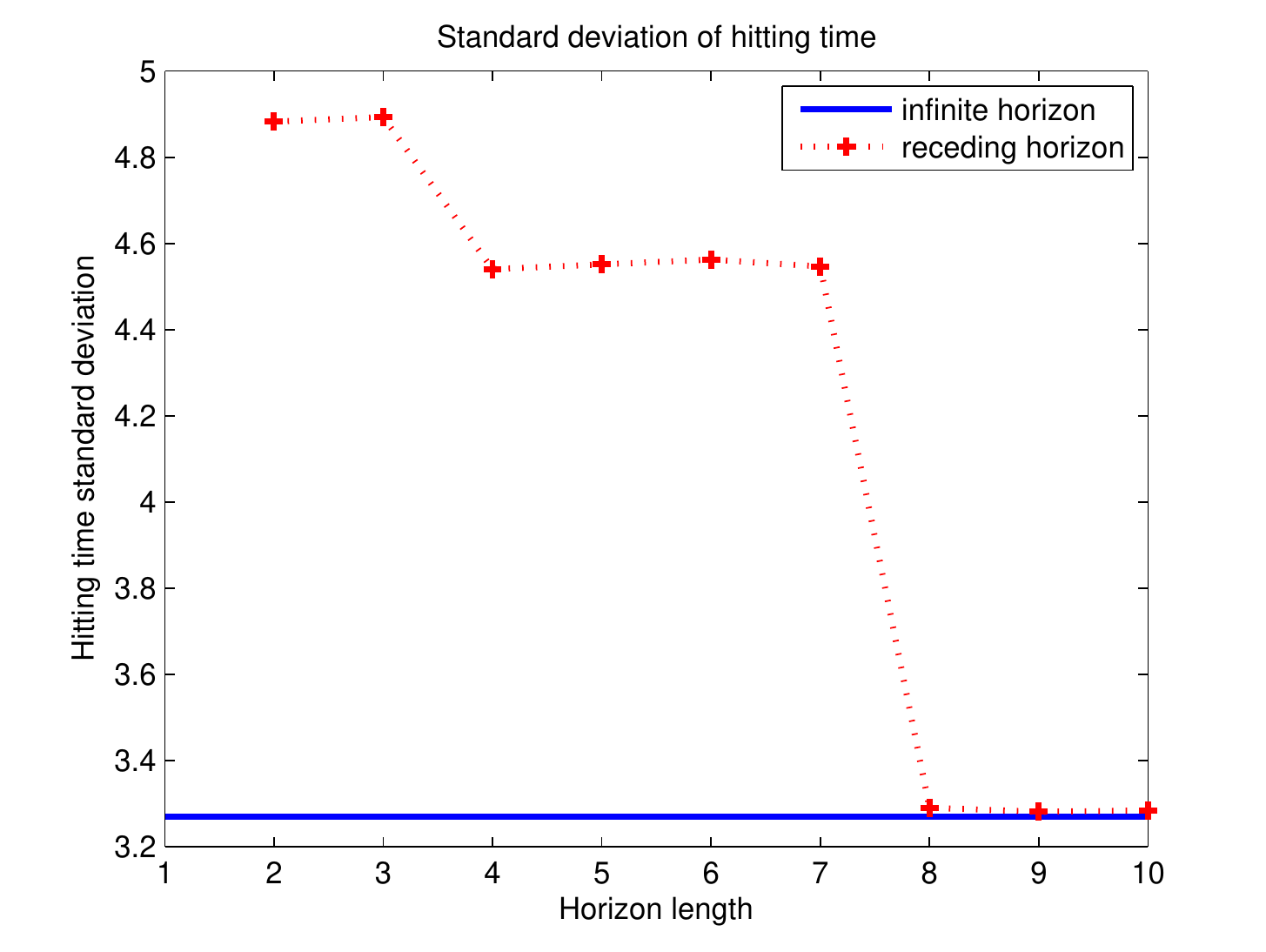}
          \caption{Hitting time average and standard-deviation}
          \label{f:time}
        \end{figure}

        Figure~\ref{f:cost} shows the average and the standard-deviation of the accumulated costs over $2\times 10^5$ Monte Carlo runs, with the initial population level at state $1$. Similarly, Figure~\ref{f:time} shows the average and the standard-deviation of the time steps needed for the recovery into the target state $4$. The results suggest that for the rolling horizon policy to match the optimal infinite horizon one, a horizon length of at least $8$ should be used. Smaller horizons provide sub-optimal policies (with respect to the infinite horizon one), with the sub-optimality gap reducing as the horizon length increases. Note that the case of $N = 1$ is not included in the data; this is because for horizon length of $1$ the optimal policy is to harvest while the system is at state $1$, leading to an $\infty$ cost and recovery time, which does not allow the system to ever recover to state $4$.

	\section{Future Work}
	\label{s:concl}
		We established in~\secref{s:EDC} that the optimal value function $V^\star$ is the minimal solution of the $\alpha$-discounted cost optimality equation~\eqref{e:alphadcoe}. However, obtaining analytical expression of the optimal value function $V^\star$ is difficult, particularly due to the integration over a subset $X\setmin K$ of the state space. Obtaining good approximations of $V^\star$ is of vital importance, and will be reported in subsequent articles.

		It is interesting to note that our basic framework of stochastic model-predictive control (described in~\secref{s:intro}) naturally leads to a partitioning of the state-space with different dynamics in each partition; thus, the controlled system may be viewed as a stochastic hybrid system. One of the basic questions in this context is that of stability of the controlled system, and in view of the fact that in general there will be infinitely many excursions of the state outside the safe set, establishing any stability property is a challenging task. Classical Lyapunov-based methods are difficult to apply directly precisely because of the infinitely many state-dependent switches between multiple regimes, each with different dynamics. However, excursion-theory of Markov processes~\cite{ref:blumenthal1992} enables us to establish certain stability properties of quite general stochastic hybrid systems with state-dependent switching; some of these results are reported in~\cite{ref:palExcur}.

	\section*{Acknowledgments}
		The authors are grateful to Vivek S. Borkar, On\'esimo Hern\'andez-Lerma, and Sean P. Meyn for illuminating discussions and pointers to relevant literature. They also thank the anonymous reviewers for their helpful comments.

%\bibliographystyle{siam}
%\bibliography{../references}

	\begin{appendix}
	\section{Proof of Theorem \ref{t:rh}}
		\begin{proof}[Proof of Theorem {\rm \ref{t:rh}}]
			For brevity of notation in this proof, we let $\wh\pi \Let \wh f_N^\infty$, and let $\wh\pi_{i:j}$ denote the (ordered) elements of the policy $\wh\pi$ from stage $i$ through $j$ for $j > i$. The first inequality in~\eqref{e:keyrh} is trivial because $v_{N+1}(x) \le V^\star(x) \le V\bigl(\wh f_N^\infty, x\bigr)$ for all $x\in X\setmin K$. Before the proof of the second inequality in~\eqref{e:keyrh}, let us fix some notation. Pick $N\in\Nz$. For $n\in\Nz$, a policy $\pi_{n:n+N}$ for stages $n$ through $n+N$, and $i\in\{n, \ldots, n+N\}$, let $\pi_{n:n+N}(i)$ denote $i$-th element of the policy $\pi_{n:n+N}$. Also, let $Q\bigl(\cdot\big| x, \pi_{n:n+N}\bigr)$ denote the sub-stochastic kernel\footnote{Recall that $Q(\cdot|\cdot)$ is a \emph{sub-stochastic kernel} on $X\setmin K$ given $Y$ if $Q(B|\cdot)$ is a measurable function on $Y$ for each $B\in\Borelsigalg{X}$, and $Q(\cdot|y)$ is a measure on $X$ with $Q(X|y) \le 1$ for each $y\in Y$.} defined for $x\in X\setmin K$ by
			\begin{align*}
				Q\bigl(B\big|x, \pi_{n:n+N}\bigr) \Let &{} \int_{X\setmin K}Q\bigl(\mrm d\xi_0\big|x, \pi_{n:n+N}(n)\bigr)\cdots\int_{X\setmin K}Q\bigl(\mrm d\xi_N\bigr|\xi_{N-1}, \pi_{n:n+N}(n+N)\bigr)\indic{B}(\xi_N)
			\end{align*}
			for $B\in\Borelsigalg{X\setmin K}$.

			Let $\pi^\star_{n:n+N}$ be an optimal policy for stages $n$ through $n+N$, i.e., let $\pi^\star_{n:n+N}$ attain the infimum in~\eqref{e:rhcf}. Fix $x\in X\setmin K$. Let $\zeta_{n+1:n+N+1}$ be an $(N+1)$-period policy starting from stage $n+1$, such that its first $N$ elements are identical to the last $N$ elements of $\pi^\star_{n:n+N}$, i.e., $\zeta_{n+1:n+N+1}(j) = \pi^\star_{n:n+N}(j)$ for $j=n+1, \ldots, n+N$. By optimality of $\pi^\star_{n:n+N}$ we have
			\begin{multline*}
				\mathsf E^{\zeta_{n+1:n+N+1}}_x\!\left[\sum_{i=n+1}^{n+N+1} \alpha^i c(x_i, a_i)\indic{\{i< \tau\}}\left.\vphantom{\sum_{i=n+1}^{n+N+1}}\right|x_{(n+1)\mn(\tau-1)}\right]\\
				\ge \mathsf E^{\pi^\star_{n+1:n+N+1}}_x\!\left[\sum_{i=n+1}^{n+N+1} \alpha^i c(x_i, a_i)\indic{\{i < \tau\}}\left.\vphantom{\sum_{i=n+1}^{n+N+1}}\right|x_{(n+1)\mn(\tau-1)}\right].
			\end{multline*}
			Since $\wh\pi_{n:n+N}(n) = \pi^\star_{n:n+N}(n)$ by construction, conditional on $x_{n\mn(\tau-1)} = x'\in X\setmin K$,
			\begin{multline}
				\label{e:keyineq}
				\int_{X\setmin K}Q\bigl(\mrm dy\big|x', \pi^\star_{n:n+N}(n)\bigr)\mathsf E^{\zeta_{n+1:n+N+1}}\!\left[\sum_{i=n+1}^{n+N+1} \alpha^i c(x_i, a_i)\indic{\{i < \tau\}}\left.\vphantom{\sum_{i=n+1}^{n+N+1}}\right|x_{(n+1)\mn(\tau-1)} = y\right]\ge\\
				\int_{X\setmin K}Q\bigl(\mrm dy\big|x', \wh\pi_{n:n+N}(n)\bigr)\mathsf E^{\pi^\star_{n+1:n+N+1}}_x\!\left[\sum_{i=n+1}^{n+N+1} \alpha^i c(x_i, a_i)\indic{\{i < \tau\}}\left.\vphantom{\sum_{i=n+1}^{n+N+1}}\right|x_{(n+1)\mn(\tau-1)} = y\right].
			\end{multline}
			By definition of $\zeta$ we have
			\begin{align*}
				\mathsf E^{\zeta_{n+1:n+N+1}}\!& \left[\sum_{i=n+1}^{n+N+1} \alpha^i c(x_i, a_i)\indic{\{i < \tau\}}\left.\vphantom{\sum_{i=n+1}^{n+N+1}}\right|x_{(n+1)\mn(\tau-1)}\right]\\
				& = \mathsf E^{\zeta_{n+1:n+N+1}}\!\left[\sum_{i=n+1}^{n+N} \alpha^i c(x_i, a_i)\indic{\{i < \tau\}}\left.\vphantom{\sum_{i=n+1}^{n+N+1}}\right|x_{(n+1)\mn(\tau-1)}\right] \\
				& \qquad+ \mathsf E^{\zeta_{n+1:n+N+1}}\!\left[\alpha^{n+N+1} c(x_{n+N+1}, a_{n+N+1})\indic{\{n+N+1 < \tau\}}\big|x_{(n+1)\mn(\tau-1)}\right],
			\end{align*}
			and the right-hand side equals%letting $\theta_{n+1:n+N}$ be an $N$-period policy starting from stage $n+1$ defined as $\theta_{n+1:n+N}(j) \Let \pi^\star_{n:n+N}(j)$ for $j=n+1, \ldots, n+N$, we rewrite the right-hand side as
			\begin{multline*}
				\mathsf E^{\pi^\star_{n:n+N}}\!\left[\sum_{i=n+1}^{n+N} \alpha^i c(x_i, a_i)\indic{\{i < \tau\}}\left.\vphantom{\sum_{i=n+1}^{n+N+1}}\right|x_{(n+1)\mn(\tau-1)}\right] \\
				+ \mathsf E^{\zeta_{n+1:n+N+1}}\!\left[\alpha^{n+N+1} c(x_{n+N+1}, a_{n+N+1})\indic{\{n+N+1 < \tau\}}\big|x_{(n+1)\mn(\tau-1)}\right].
			\end{multline*}
			In conjunction with~\eqref{e:keyineq} and conditional on $x_{n\mn(\tau-1)} = x'\in X\setmin K$, we have
			\begin{multline*}
				\int_{X\setmin K}Q\bigl(\mrm dy\big|x', \pi^\star_{n:n+N}(n)\bigr)\mathsf E^{\pi^\star_{n:n+N}}\!\left[\sum_{i=n+1}^{n+N} \alpha^i c(x_i, a_i)\indic{\{i < \tau\}}\left.\vphantom{\sum_{i}^N}\right|x_{(n+1)\mn(\tau-1)}=y\right]\\
					+ \int_{X\setmin K} Q\bigl(\mrm dy\big|x', \pi^\star_{n:n+N}(n)\bigr)\mathsf E^{\zeta_{n+1:n+N+1}}\bigl[\alpha^{n+N+1} c(x_{n+N+1}, a_{n+N+1})\indic{\{n+N+1 < \tau\}}\big|x_{(n+1)\mn(\tau-1)}=y\bigr]\\
					\ge \int_{X\setmin K} Q\bigl(\mrm dy\big|x', \wh\pi_{n:n+N}(n)\bigr)\mathsf E^{\pi^\star_{n+1:n+N+1}}\!\left[\sum_{i=n+1}^{n+N+1}\alpha^i c(x_i, a_i)\indic{\{i < \tau\}}\left.\vphantom{\sum_i^N}\right|x_{(n+1)\mn(\tau-1)}=y\right].
			\end{multline*}
			To wit, conditional on $x_{n\mn(\tau-1)} = x'\in X\setmin K$,
			\begin{multline*}
				\mathsf E^{\pi^\star_{n:n+N}}\!\left[\sum_{i=n}^{n+N} \alpha^i c(x_i, a_i)\indic{\{i < \tau\}}\left.\vphantom{\sum_{i}^N}\right|x_{n\mn(\tau-1)} = x'\right] - \mathsf E^{\pi^\star_{n:n+N}}\!\left[\alpha^n c(x_n, a_n)\indic{\{i < \tau\}}\big|x_{n\mn(\tau-1)} = x'\right]\\
					+ \alpha^{n+N+1} \int_{X\setmin K} \!\!Q\bigl(\mrm dy\big|x', \pi^\star_{n:n+N}(n)\bigr)\mathsf E^{\zeta_{n+1:n+N+1}}\bigl[c(x_{n+N+1}, a_{n+N+1})\indic{\{n+N+1 < \tau\}}\big|x_{(n+1)\mn(\tau-1)}=y\bigr]\\
					\ge \int_{X\setmin K} Q\bigl(\mrm dy\big|x', \wh\pi_{n:n+N}(n)\bigr)\mathsf E^{\pi^\star_{n+1:n+N+1}}\!\left[\sum_{i=n+1}^{n+N+1}\alpha^i c(x_i, a_i)\indic{\{i < \tau\}}\left.\vphantom{\sum_i^N}\right|x_{(n+1)\mn(\tau-1)}=y\right].
			\end{multline*}
			Let $\zeta_{n+1:n+N+1}(n+N+1)(\cdot)$ be a selector that attains the minimal value of
			\[
				\alpha^{n+N+1}\int_{X\setmin K}Q\bigl(\mrm dy\big|x', \pi^\star_{n:n+N}(n)\bigr)\mathsf E^{\zeta_{n+1:n+N+1}}\bigl[c(x_{n+N+1}, a_{n+N+1})\indic{\{n+N+1 < \tau\}}\big|x_{(n+1)\mn(\tau-1)}=y\bigr]
			\]
			whenever $x'\in X\setmin K$, and let the corresponding minimal value be denoted by $e_n(x')$; clearly $e_n$ is well-defined on $X\setmin K$, and is a measurable function of $x'$. With this notation, the last inequality becomes
			\begin{multline*}
				\mathsf E^{\pi^\star_{n:n+N}}\!\left[\sum_{i=n}^{n+N}\alpha^i c(x_i, a_i)\indic{\{i < \tau\}}\left.\vphantom{\sum_i^N}\right|x_{n\mn(\tau-1)} = x'\right] - \mathsf E^{\pi^\star_{n:n+N}}\bigl[\alpha^n c(x_n, a_n)\indic{\{n < \tau\}}\big|x_{n\mn(\tau-1)} = x'\bigr]\\
				+ e_n(x') \ge \int_{X\setmin K} Q\bigl(\mrm dy\big|x', \wh\pi_{n:n+N}(n)\bigr) \mathsf E^{\pi^\star_{n+1:n+N+1}}\!\left[\sum_{i=n+1}^{n+N+1}\alpha^i c(x_i, a_i)\indic{\{i < \tau\}}\left.\vphantom{\sum_i^N}\right|x_{(n+1)\mn(\tau-1)}=y\right]
			\end{multline*}
			whenever $x'\in X\setmin K$. Therefore,
			\begin{multline*}
				\int_{X\setmin K}Q\bigl(\mrm dy\big|x, \wh\pi_{0:n-1}\bigr) \mathsf E^{\pi^\star_{n:n+N}}\!\left[\sum_{i=n}^{n+N}\alpha^i c(x_i, a_i)\indic{\{i < \tau\}}\left.\vphantom{\sum_i^N}\right|x_{n\mn(\tau-1)}=y\right]\\
					- \int_{X\setmin K}Q\bigl(\mrm dy\big|x, \wh\pi_{0:n-1}\bigr) \mathsf E^{\pi^\star_{n:n+N}}\bigl[\alpha^n c(x_n, a_n)\indic{\{n < \tau\}}\big|x_{n\mn(\tau-1)}\bigr] + \int_{X\setmin K}Q\bigl(\mrm dy\big|x, \wh\pi_{0:n-1}\bigr) e_n(y)\\
					\ge \int_{X\setmin K}Q\bigl(\mrm dy\big|x, \wh\pi_{0:n}\bigr) \mathsf E^{\pi^\star_{n+1:n+N+1}}\!\left[\sum_{i=n+1}^{n+N+1}\alpha^i c(x_i, a_i)\indic{\{i < \tau\}}\left.\vphantom{\sum_i^N}\right|x_{(n+1)\mn(\tau-1)}=y\right].
			\end{multline*}
			Rearranging and summing over $n$ we arrive at
			\begin{multline}
			\label{e:ineq1}
				\sum_{n=0}^\infty \int_{X\setmin K}Q\bigl(\mrm dy\big|x, \wh\pi_{0:n-1}\bigr) \mathsf E^{\pi^\star_{n:n+N}}\!\left[\alpha^n c(x_n, a_n)\indic{\{n < \tau\}}\left.\vphantom{\sum}\right|x_{n\mn(\tau-1)}=y\right]\\
				\le \sum_{n=0}^\infty\left(\alpha^n\int_{X\setmin K}Q\bigl(\mrm dy\big|x, \wh\pi_{0:n-1}\bigr)\mathsf E^{\pi^\star_{n:n+N}}\!\left[\sum_{i=n}^{n+N}\alpha^{i-n} c(x_i, a_i)\indic{\{i < \tau\}}\left.\vphantom{\sum_i^N}\right| x_{n\mn(\tau-1)}=y\right]\right.\\
				- \left.\alpha^{n+1}\int_{X\setmin K} Q\bigl(\mrm dy\big|x, \wh\pi_{0:n}\bigr)\mathsf E^{\pi^\star_{n+1:n+N+1}}\!\left[\sum_{i=n+1}^{n+N+1}\alpha^{i-n-1} c(x_i, a_i)\indic{\{i < \tau\}}\left.\vphantom{\sum_i^N}\right| x_{(n+1)\mn(\tau-1)}=y\right]\right)\\
				+ \sum_{n=0}^\infty \int_{X\setmin K} Q\bigl(\mrm dy\big|x, \wh\pi_{0:n-1}\bigr) e_n(y).
			\end{multline}
			In~\eqref{e:ineq1} we have employed the notation $\int_{X\setmin K}Q\bigl(\mrm dy|x, \pi_{0:-1}\bigr) g(y) \Let g(x)$ for any policy $\pi$. We observe that the left-hand side of~\eqref{e:ineq1} is just $\mathsf E^{\wh\pi}_x\!\left[\sum_{i=0}^{\tau-1} \alpha^i c(x_i, a_i)\right]$. By Assumption~\ref{a:further}(i),
			\[
				\mathsf E^{\pi^\star_{n:n+N}}\!\left[\sum_{i=n}^{n+N}\alpha^{i-n} c(x_i, a_i)\indic{\{i < \tau\}}\left.\vphantom{\sum_i^N}\right| x_{n\mn(\tau-1)}=y\right] \le \ol c\mathsf E^{\pi^\star_{n:n+N}}\!\left[\sum_{i=n}^{n+N}\alpha^{i-n} w(x_i)\indic{\{i < \tau\}}\left.\vphantom{\sum_i^N}\right| x_{n\mn(\tau-1)}=y\right],
			\]
			and by Assumption~\ref{a:further}(ii),
			\[
				\mathsf E^{\pi^\star_{n:n+N}}\!\left[\sum_{i=n}^{n+N}\alpha^{i-n} w(x_i)\indic{\{i < \tau\}}\left.\vphantom{\sum_i^N}\right| x_{n\mn(\tau-1)}=y\right] \le w(y)\sum_{i=n}^{n+N}\gamma^{i-n}.
			\]
			We notice that since $c\ge 0$, the first series on the right-hand side of~\eqref{e:ineq1} is at most
			\begin{equation}
			\label{e:ineq2}
			\begin{aligned}
				\ol c \sum_{n=0}^\infty \alpha^n\sum_{i=n}^{n+N}\gamma^{i-n} \int_{X\setmin K}Q\bigl(\mrm dy\big|x, \wh\pi_{0:n-1}\bigr) w(y).
			\end{aligned}
			\end{equation}
			For a fixed $n\in\Nz$, the quantity $\int_{X\setmin K}Q\bigl(\mrm dy|x, \wh\pi_{0:n}\bigr)w(y)$ is at most $\beta^{n+1} w(x)$ in view of Assumption~\ref{a:further}(ii) and the definition of the stochastic kernel $Q\bigl(\cdot\big|x, \pi_{n:n+N}\bigr)$ at the beginning of this proof. Therefore,
			\begin{align*}
				\sum_{n=0}^\infty & \alpha^n\sum_{i=n}^{n+N}\gamma^{i-n} \int_{X\setmin K}Q\bigl(\mrm dy\big|x, \wh\pi_{0:n-1}\bigr) w(y) \le \sum_{n=0}^\infty\alpha^n\sum_{i=n}^{n+N}\gamma^{i-n}\beta^n w(x)\\
				& \le w(x) \left(\frac{1-\gamma^{N+1}}{1-\gamma}\right) < \infty.
			\end{align*}
			This shows that series in~\eqref{e:ineq2} is summable. Hence, cancellations of the telescopic terms in the first series on the right-hand side of~\eqref{e:ineq1} are justified. The inequality in~\eqref{e:ineq1} now simplifies to
			\begin{multline}
			\label{e:ineq3}
				\mathsf E^{\wh\pi}_x\!\left[\sum_{i=0}^{\tau-1} \alpha^i c(x_i, a_i)\right] \le \mathsf E^{\pi^\star_{0:N}}_x\!\left[\sum_{i=0}^{(N+1)\mn(\tau-1)}\!\!\alpha^{i} c(x_i, a_i)\right] + \sum_{n=0}^\infty \int_{X\setmin K} Q\bigl(\mrm dy\big|x, \wh\pi_{0:n-1}\bigr) e_n(y).
			\end{multline}
			By Assumption~\ref{a:further}(ii) and the definition of $e_n$, conditional on $x_{n\mn\tau} = x'\in X\setmin K$,
			\begin{align*}
				& e_n(x')\\
				& \le \alpha^{n+N+1}\int_{X\setmin K} \!\!Q\bigl(\mrm dy\big| x', \pi^\star_{n:n+N}(n)\bigr)\mathsf E^{\zeta_{n+1:n+N+1}}\bigl[c(x_{n+N+1}, a_{n+N+1})\indic{\{n+N+1 < \tau\}}\big|x_{(n+1)\mn(\tau-1)} = y\bigr]\\
				& \le \ol c w(x')\alpha^n \gamma^{N+1}.
			\end{align*}
			Substituting the last inequality in~\eqref{e:ineq3} we arrive at
			\[
				\mathsf E^{\wh\pi}_x\!\left[\sum_{i=0}^{\tau-1} \alpha^i c(x_i, a_i)\right] \le \mathsf E^{\pi^\star_{0:N}}_x\!\left[\sum_{i=0}^{(N+1)\mn(\tau-1)}\!\!\alpha^{i} c(x_i, a_i)\right] + \frac{\ol c\gamma^{N+1}}{1-\gamma}w(x),
			\]
			which is the second bound in~\eqref{e:keyrh}. The inequality~\eqref{e:sloppyrh} follows immediately from the fact that $V^\star \ge v_n$ for every $n\in\N$.
		\end{proof}
	\end{appendix}

\def\cprime{$'$}

%\bibliographystyle{siam}
%\bibliography{../references}

\end{document}